\documentclass[11pt]{article}
\usepackage{amssymb}
\usepackage{amsmath}
\usepackage{theorem}
\usepackage{algorithm}
\usepackage{graphicx}
\usepackage[usenames]{color}
\usepackage{multirow}

\newcommand{\R}{\mathbb{R}}
\newcommand{\N}{\mathbb{N}}

\newcommand{\ve}[1]{ #1}
\newtheorem{Thm}{Theorem}[section]
\newtheorem{Def}{Definition}[section]
\newtheorem{Cor}{Corollary}[section]
\newtheorem{Lemma}{Lemma}[section]

  {\end{list}}

\theoremheaderfont{\scshape}
\makeatletter

  \gdef\listctr{list\romannumeral\the\@listdepth}\expandafter
  
\makeatother

\newenvironment{AlgorithmSteps}[1][1]{%
  \begin{list}{\csname label\listctr\endcsname}{%
      \usecounter{\listctr}
      
      \settowidth{\labelwidth}{\textsc{Step\ #1.}}%
      \setlength{\leftmargin}{\labelwidth}\addtolength{\leftmargin}{\labelsep}}}%
  {\end{list}}

\def\proof{{\it Proof. }}
\def\x{ {\ve x}}

\def\y{ {\ve y}} 
\def\u{ {\ve u}}
\def\w{{\ve w}}
\def\z{ {\ve z}}

\def\p{{\ve p}}
\def\xk{ \x^{(k)}}
\def\yk{ \y^{(k)}}

\def\uk{ \u^{(k)}}
\def\zk{ \z^{(k)}}\def\zkk{ \z^{(k+1)}}
\def\xs{\x^*}
\def\yk{ \y^{(k)}}

\def\xkk{ \x^{(k+1)}}
\def\xkm{ \x^{(k-1)}}

\def\ukk{ \u^{(k+1)}}

\def\dom{\mbox{dom}}
\def\R{\mathbb R}

\def\endproof{\hfill$\square$\vspace{0.3cm}\\}
\def\prox{{\mbox{prox}}}
\newcommand\D{{\mathcal D}}
\newcommand{\stl}{\alpha}
 \newcommand{\projj}[1]{\p_{\stl,D}({#1})}
\newcommand{\KL}{\operatorname{KL}}
\newcommand{\HS}{\operatorname{HS}}
\newcommand{\projjk}[2][k]{\p_{\stl_{#1},D_{#1}}({#2})}

\newcommand{\xkb}{\x^{(k-1)}}
\newcommand{\tk}{t_k}\newcommand{\tkm}{t_{k-1}} \newcommand{\ak}{\stl_k}\newcommand{\akk}{\stl_{k+1}}
\newcommand{\vk}{v_k}\newcommand{\vkk}{v_{k+1}}\newcommand{\Dk}{{D_k}}
 
\newcommand{\thk}{\theta_k}\newcommand{\thkb}{\theta_{k-1}}
\newcommand{\sk}{s_k}\newcommand{\skk}{s_{k+1}}
\def\x{ {x}}\def\y{ {y}}
\def\xkk{ {\x^{(k+1)}}} 
\def\xk{ \x^{(k)} }\def\yk{ \y^{(k)} }
\def\xkm{ \x^{(k-1)} }
\def\prox{\mathrm{prox}} 
\newcommand{\argmin}{\operatornamewithlimits{argmin}}
\def\dom{\mbox{dom}}

\begin{document}
\title{A variable metric forward--backward method with extrapolation\thanks{This work has been partially supported by MIUR under the project FIRB - Futuro in Ricerca 2012, contract RBFR12M3AC. The Italian GNCS - INdAM is also acknowledged.}}
\author{S. Bonettini, F. Porta, V. Ruggiero}
\maketitle
\begin{abstract}
Forward-backward methods are a very useful tool for the minimization of a functional
given by the sum of a differentiable term and a nondifferentiable one and their investigation has experienced several efforts from many researchers in the last decade. In this paper we focus on the convex case and, inspired by recent approaches for accelerating first-order iterative schemes,
we develop a scaled inertial forward-backward algorithm which is based on a metric changing at each iteration and on a suitable extrapolation step.
Unlike standard forward-backward methods with extrapolation, our scheme is able to handle functions whose domain is not the entire space. Both {an ${\mathcal O}(1/k^2)$ convergence rate estimate on the objective function values and the convergence of the sequence of the iterates} are proved. Numerical experiments on several {test problems arising from image processing, compressed sensing and statistical inference}
show the {effectiveness} of the proposed method in
comparison to well performing {state-of-the-art} algorithms.
\end{abstract}
\section{Introduction}\label{sec:1}
In this paper we are interested in solving the optimization problem
\begin{equation}\label{minf}
\min_{\x\in\R^n} F(\x) \equiv f(\x) + g(\x)
\end{equation}
where $f$ and $g$ are proper, convex and lower semicontinuous functions from $\R^n$ to $\R\cup \{\infty\}$. Moreover, we assume that $f$ is differentiable with Lipschitz continuous gradient on a suitable closed, convex set $Y\subseteq \dom(f) = \{\x\in\R^n:f(\x)<\infty\}$, such that
$$
\dom(f) \supseteq Y\supseteq \dom(g)
$$
We also suppose that $g$ is bounded from below over its domain and problem \eqref{minf} admits at least a solution.
Formulation \eqref{minf} includes also constrained problems over a closed convex set $\Omega\subset \R^n$ where $f$ has Lipschitz continuous gradient: indeed, the constraints defined by $\Omega$ can be inserted into the model by adding to $g$ the indicator function of the feasible set itself, i.e.
$$
\min_{x\in \Omega} f(x)+g(x) = \min_{x\in \R^n} f(x)+g(x)+\iota_{\Omega}(x) 
$$
where
$$
\iota_{\Omega}(x) =
\left\{
\begin{array}{cc}
 0 & {\rm if} \ \ x\in\Omega\\
 +\infty & {\rm otherwise}
\end{array}
\right.
$$
Problem \eqref{minf} is relevant in various domains of applied science such as signal and image processing, statistical inference and machine learning. A typical feature of these applications is the large number of variables, which makes the class of first order methods very attractive.
In this class, forward-backward methods \cite{Combettes-Pesquet-2009,Combettes-Wajs-2005} are especially suited for problem \eqref{minf}, since they exploit the decomposition of the objective function in a differentiable term and a nondifferentiable one. The general forward-backward iteration is given by
\begin{equation}\nonumber
\xkk = \xk + \lambda_k(\prox_{\alpha_k g}(\xk-\alpha_k \nabla f(\xk)) - \xk)  \, ,
\end{equation}
where $\lambda_k$, $\alpha_k$ are positive parameters controlling the steplength and $\prox_{\phi}(\cdot)$ is the proximity operator associated to the convex function $\phi${, 
defined as}
\begin{equation}\nonumber
\prox_\phi(y) =  \argmin_{\x\in\R^n} \ \phi(\x) + \frac 1 2 \|\x-\y\|^2
\end{equation}
Forward-backward methods are easy to implement and have well studied convergence properties.
On the other hand, it is well known that they can exhibit a poor convergence rate, especially when a high accuracy is required.

In the recent literature, we can find two different approaches aiming to improve the convergence speed of forward-backward methods.  They are both described below.
\paragraph{Inertial/estrapolation techniques.} This approach consists in adding an extrapolation step to the basic forward-backward iteration, yielding a multistep algorithm, called also heavy ball or inertial method \cite[p.65]{Polyak-1987}. The idea of inertial methods became very popular in the last decade, in view of Nesterov's work \cite{Nesterov-2005} and it has been further developed in \cite{Beck-Teboulle-2009b}, where the authors propose the following variant
    \begin{eqnarray}
    \yk  &=& \xk+ \beta_k(\xk-\xkm)\label{FISTA1}\\
    \xkk &=& \prox_{\alpha_k g}(\yk-\alpha_k\nabla f(\yk))\label{FISTA2} \, .
    \end{eqnarray}
In \cite{Beck-Teboulle-2009b,Bertsekas-2012}, the convergence of method \eqref{FISTA1}--\eqref{FISTA2} is investigated by showing that for suitable sequences of parameters $\{\alpha_k\}$ and $\{\beta_k\}$ (with $\lim_k \beta_k = 1$) one has
$F(\xk) - F^* = {\mathcal O}\left(\frac 1 {k^2}\right)$,
where $F^*$ is the optimal value of the objective function. Recently, under additional assumptions on the sequences $\{\alpha_k\}$ and $\{\beta_k\}$, Chambolle and Dossal in \cite{Chambolle-Dossal-2014} proved the convergence of the iterates $\{\xk\}$, while in \cite{Villa-etal-2013} the authors propose a variant of the inertial scheme where the proximal point \eqref{FISTA2} can be computed inexactly. We also mention the recent work \cite{Ochs-etal-2014} where inertial forward-backward algorithms are analyzed when $f(\x)$ is not convex.

A drawback in the use of method \eqref{FISTA1}--\eqref{FISTA2} is that it may be unfeasible when $\dom(f)$ does not coincide with the whole space $\R^n$, since the point $\yk$ computed in \eqref{FISTA1} does not necessarily belong to $\dom(f)$.
\paragraph{Variable metric/scaling techniques.} In a variable metric forward-backward algorithm, the underlaying metric may change at each iteration by means of suitable symmetric positive definite scaling matrices multiplying the gradient of $f$  and also involved in the definition of the proximity operator. The expected advantage in using a variable metric is an improved capability to capture the local features of problem \eqref{minf}, possibly leading to an improvement of the convergence speed (think for example to the Newton's method). In \cite{Combettes-Vu-2013,Combettes-Vu-2014}, the authors propose and analyze the following variable metric forward-backward algorithm
\begin{equation}\label{SFB}
\xkk = \xk + \lambda_k(\prox_{\alpha_k g}^{D_k}(\xk-\alpha_k D_k^{-1}\nabla f(\xk)) - \xk)
\end{equation}
where $\{D_k\}$ is a user supplied sequence of symmetric positive definite matrices and $\prox_{\alpha_k g}^{D_k}(\y)$ is defined as
\begin{equation}\label{scal-prox}
\prox_{\alpha_k g}^{D_k}(y) =  \argmin_{\x\in\R^n}\ g(\x) + \frac 1 {2\alpha_k} (\x-\y)^TD_k(\x-\y)
\end{equation}
The authors also devise conditions on the sequence $\{D_k\}$ ensuring the convergence of $\{\xk\}$. Method \eqref{SFB}, equipped by an Armijo line--search for the computation of $\lambda_k$, has been extensively studied in the papers \cite{Birgin-etal-2003, Bonettini-Prato-2015a, Bonettini-etal-2009} {for constrained minimization,} when $g(\x)$ reduces to the indicator function of a closed convex subset of $\R^n$. The convergence rate on the objective function values in this case is only linear, i.e.
$F(\xk) - F^* = {\mathcal O}\left(\frac 1 {k}\right)$
(see \cite{Beck-Teboulle-2009b,Bonettini-Prato-2015a}). However, in spite of the theoretical convergence rate, a suitable combination of the stepsize parameter $\alpha_k$ and the scaling matrix $D_k$ can allow method \eqref{SFB} to reach practical performances which are comparable with those of \eqref{FISTA1}--\eqref{FISTA2} \cite{Bonettini-Prato-2015a,Porta-Prato-Zanni-2015}.
\paragraph{Main contribution.} In this paper we propose an algorithm combining the two acceleration techniques described above, designing an original variable metric forward-backward method with extrapolation.

In particular, we address the case where $\dom(f)$ is a proper subset of $\R^n$ and we devise suitable conditions on the stepsize parameters and on the scaling matrices sequence to ensure both the convergence of the iterates sequence $\{\xk\}$ and the ${\mathcal O}\left(\frac 1 {k^2}\right)$ rate {for} the objective function values.

The effectiveness of the proposed method is evaluated by means of a comparison with other state-of-the-art algorithms, on several optimization problems of the form \eqref{minf}, arising
from different real-life applications such as image deblurring, compressed sensing and probability density estimation.

The plan of the paper is the following. In Section \ref{sec:2} we collect some definitions and introductory results. Section \ref{sec:algodescription} is devoted to the description of the proposed algorithm while the convergence rate analysis is performed in Section \ref{sec:rate} and the convergence of the iterates to a minimizer of the optimization problem is proved in Section \ref{sec:conv_iter}. {This last section is strongly inspired from a recent paper by Chambolle and Dossal \cite{Chambolle-Dossal-2014}.}
Section \ref{sec:num} deals with the results of the numerical experiments we performed on {some test problems} arising in image and signal processing and statistical inference. 
Finally, our conclusions are given in Section \ref{sec:concl}.

\section{Notation, definitions and basic results}\label{sec:2}
We denote by $\|\cdot\|$ the Euclidean norm of a vector while $\|\cdot\|_D$ indicates the norm induced by the symmetric positive definite matrix $D$, i.e. $\|x\|_D^2=x^T Dx$. Furthermore,
in the subspace $\mathcal{S}_n(\R)$ of the symmetric {real} matrices of order $n$, we consider the following Loewner partial ordering
$$\forall D_1, D_2 \in \mathcal{S}_n(\R) \quad D_1\succeq D_2 \Leftrightarrow  x^T D_1 x\geq x^T D_2 x \quad \forall  x\in \R^n$$
For any $\eta\in\R$, $\eta>0$ we define the set $\D_\eta\subset \mathcal{S}_n(\R)$ as the set of all positive definite matrices $D$ such that
$D \succeq \eta I$.
Clearly, if $D\in\D_\eta$, the eigenvalues of $D$ are lower bounded by $\eta$ and
for each $\u\in\R^n$, the following inequality holds
\begin{equation}\label{ine_norm}
\eta \| \u \|^2 \leq \u^T D \u = \|\u\|_D^2
\end{equation}
The following lemma states a well known property of the projection operator, whose proof runs exactly as in \cite[p.48]{Hiriart-Lemarechal-2001}
\begin{Lemma} \label{lemma:nonexpansive}
Let $\Omega\subseteq \R^n$ be a closed convex set and define the scaled Euclidean projection operator associated to $D\in \D_\eta$ as
\begin{equation}\label{scaled_euclidean_proj}
P_{\Omega,D}(\x) = \argmin_{\y\in \Omega}\|\y-\x\|^2_D
\end{equation}
for any $\x\in \R^n$.
Then, the operator \eqref{scaled_euclidean_proj} is nonexpansive with respect to the norm induced by the matrix $D$, i.e.
\begin{equation}\label{nonexpansive}
\|P_{\Omega,D}(\x) - P_{\Omega,D}(\z)\|_D\leq \|\x-\z\|_D
\end{equation}
for all $\x,\z\in \R^n$.
\end{Lemma}
%
For every $\x\in Y$ and $\y\in \R^n$ we define
\begin{equation}\label{def:ell}
\ell(\y;\x) = f(\x)+\nabla f(\x)^T(\y-\x)
\end{equation}
and
\begin{equation}\label{def:q}
q(\y;\x) = \ell(\y;\x) + g(\y)
\end{equation}
\begin{Def}
A smooth function $f:\R^n\to \R\cup\{\infty\}$  has $L$-Lipschitz continuous gradient on the set $\Omega\subseteq \dom(f)$ if there exists $L>0$ such that
\begin{equation}\nonumber
\|\nabla f(\x) -\nabla f(\y)\| \leq L \|\x-\y\|, \ \ \ \forall \x,\y\in \Omega.
\end{equation}
\end{Def}
We recall also the following result for smooth functions with $L$-Lipschitz continuous gradient. 
\begin{Lemma}\label{cor:1}
Let $f:\R^n\rightarrow \R\cup\{\infty\}$ be a continuously differentiable function with {$L$-}Lipschitz continuous gradient on $Y\subseteq \R^n$ 
and $D\in\D_\eta$. Then, for every $\x,\y\in Y$ we have
\begin{equation}\label{Lipschitz}
f(\y)\leq \ell(\y;\x) +\frac 1 {2\alpha} \|\x-\y\|^2_D
\end{equation}
for all $\alpha\leq \eta/L$.
\end{Lemma}
\proof From \eqref{ine_norm} we obtain
\begin{equation}\nonumber
\ell(\y;\x) +\frac 1 {2\alpha} \|\x-\y\|^2_D\geq \ell(\y;\x) +\frac{\eta}{2\alpha} \|\x-\y\|^2\geq f(\y)
\end{equation}
where the rightmost inequality holds for $\frac{\alpha}{\eta}\leq1/L$, thanks to Lemma 6.9.1 in \cite{Bertsekas-2012} (see also \cite[Lemma 2.1]{Beck-Teboulle-2009b}).\endproof
\begin{Def}
Let $g:\R^n\to \R\cup\{\infty\}$ be a convex function. Then, the subdifferential of $g$ at $\x\in\R^n$ is the set
\begin{equation}\nonumber
\partial g(\x) = \{\w\in \R^n: g(\y)\geq g(\x) + (\y-\x)^T\w, \ \ \forall \y\in \R^n\}
\end{equation}
A point $\x$ is a minimizer of $g$ if and only if $0\in \partial g(x)$.
\end{Def}
As a consequence of the previous definition, we have that $\x\in\R^n$ is a solution of problem \eqref{minf} if and only if $-\nabla f(\x)\in \partial g(\x)$.

Given a positive number $\alpha$ and a matrix $D\in\D_\eta$, we now define the following function
\begin{equation}\label{def:Q}
Q_{\alpha,D}(\y;\x) = q(\y;\x) + \frac{1}{2\alpha}\|\x-\y\|_D^2
\end{equation}
for $\y\in \R^n$, $\x\in Y$. The function $Q_{\alpha,D}(\cdot;\x)$ admits a unique minimizer, which will be denoted by
\begin{equation}\label{proj}
\projj\x = \argmin_{\y\in\R^n} Q_{\alpha,D}(\y;\x)
\end{equation}
Clearly, the point $\projj\x$ {belongs} to $\dom(g)$  {and, when $\projj\x=\x$, $\x$ is a minimizer of $F$}. A simple computation shows that
\begin{equation}\label{scal-proj}
\projj\x=\argmin_{\y\in\R^n} g(\y) + \frac{1}{2\stl} \left\| \y-\x + \stl D^{-1}\nabla f(\x)\right\|^2_D
\end{equation}
where it is more evident that {the introduction of} a matrix $D$ in \eqref{def:Q} induces a scaling of $\nabla f(\y)$ by $D^{-1}$. Clearly, according to \eqref{scal-prox}, we have the equivalence $\projj\x = \prox_{\alpha g}^{D}(\x-\alpha D^{-1}\nabla f(\x))$.
In the next {lemmas,} two useful properties of the operator \eqref{proj} are proved.
\begin{Lemma}\label{lemma:2}
Let  {$f:Y\rightarrow \R$} be a continuously differentiable function with $L$-Lipschitz continuous gradient, $D\in\D_\eta$ and $g$ be a convex function. Let  {$\x\in Y$} and $\y=\projj\x$. Then, {for any $z\in \R^n$},  we have
\begin{equation}\label{lemma:2b}
q(\y;\x) + \frac 1{2\alpha}\|\x-\y\|_D^2 \leq q(\z;\x) + \frac 1{2\alpha}\|\z-\x\|_D^2- \frac 1{2\alpha}\|\z-\y\|_D^2
\end{equation}
\end{Lemma}
\proof 
From the optimality conditions of the problem \eqref{proj}, we have that there exists a vector $\w\in \partial g(\y)$ such that
\begin{equation}\label{opt_cond_proj}
\nabla f(\x) + \frac 1 \alpha D(\y-\x) + w = 0
\end{equation}
Since $\w\in \partial g(\y)$, for all $\z\in \R^n$ we have that $g(\z)-g(\y)\geq \w^T(\z-\y)$, which, together with \eqref{opt_cond_proj}, implies
\begin{equation}\label{ine0}
g(\z)-g(\y) \geq \left(\nabla f(\x) + \frac 1 \alpha D(\y-\x)\right)^T(\y-\z)
\end{equation}
From {the} definition of $\ell(\y;\x)$ in \eqref{def:ell}, we can write
\begin{eqnarray*}
\ell(\y;\x) + \frac 1 {2\alpha} \|\y-\x\|^2_D&=& f(\x) + \nabla f(\x)^T(\y-\z+\z-\x)+\frac 1 {2\alpha} \|\y-\z+\z-\x\|^2_D\\
&=&\ell(\z;\x) + \frac 1 {2\alpha} \|\z-\x\|^2_D + \frac 1 {2\alpha} \|\y-\z\|^2_D  +\\
& & + \nabla f(\x)^T(\y-\z) + \frac 1 \alpha (\z-\x)^TD(\y-\z) \\
&=&\ell(\z;\x) + \frac 1 {2\alpha} \|\z-\x\|^2_D- \frac 1 {2\alpha} \|\y-\z\|^2_D +  \\
& & + \left(\nabla f(\x) + \frac 1 \alpha D(\y-\x)\right)^T(\y-\z)\\
&\leq& \ell(\z;\x) + \frac 1 {2\alpha} \|\z-\x\|^2_D- \frac 1 {2\alpha} \|\y-\z\|^2_D +g(\z)-g(\y)
\end{eqnarray*}
where the third equality is obtained by adding and subtracting $\frac{1}{\alpha} \y^T D(\y-\z)$ and the final inequality is a consequence of \eqref{ine0}. Finally, inequality \eqref{lemma:2b} follows by rearranging terms and recalling the definition \eqref{def:q}.\endproof
A direct consequence of the previous lemma is the following result.
\begin{Lemma}\label{lemma:2c}
Let $f:Y\rightarrow \R$ be a {convex}, continuously differentiable function with $L$-Lipschitz continuous gradient and $g$ be a convex function. Let  $F(\x)\equiv f(\x)+g(\x)$, $D\in \D_\eta$ and $\y=\projj\x$, $\x\in Y$. If $\alpha$ is such that $\y$ satisfies the condition \eqref{Lipschitz}, then,
for any {$z\in \dom(f)$},  we have
\begin{equation}\label{lemma:2d}
F(\y) + \frac 1{2\alpha}\|\z-\y\|_D^2 \leq F(\z) + \frac 1{2\alpha}\|\z-\x\|_D^2
\end{equation}
\end{Lemma}
\proof
Since $f$ is convex, the following inequality holds:
\begin{equation}\nonumber
F(z) \geq q(\z;\x) \ \ \forall \x,\z
\end{equation}
Therefore, from Lemma \ref{lemma:2} and from \eqref{Lipschitz} we have the result.
\endproof
\section{A {variable metric inertial forward-backward} method with backtracking}\label{sec:algodescription}
In this section we describe and analyze the proposed method, whose generic scheme is detailed in Algorithm \ref{SGEM}.
\begin{algorithm}[h]
\caption{Scaled  {inertial forward-backward} method 
with backtracking}\label{SGEM}
{Choose $\stl_0> 0$, $\delta <1 $, $x^{(0)}\in Y$. Set $\x^{(-1)}=\x^{(0)}$ and define a sequence of nonnegative numbers $\{\beta_k\}$ and a sequence of matrices $\{D_k\}$, with $D_k\in \D_\eta$, such that $\gamma=\mbox{sup}_{k\in \N} \|D_k\|<\infty$.
\begin{itemize}
\item[]\textsc{For} $k=0,1,2,...$
\begin{AlgorithmSteps}[4]
\item[1] Extrapolation:
$\yk = P_{Y,D_k}(\xk+\beta_k(\xk-\xkb))$
\item[2] Set $\stl_k = \stl_{k-1}$, $i_k=0$
\item[3] Set $\xk_+ = \projjk\yk$
\item[4] If
\begin{equation}\nonumber
f(\xk_+) \leq \ell(\xk_+;\yk) + \frac{1}{2\ak}\|\yk-\xk_+\|^2_\Dk
\end{equation}
go to Step 5.\\
else set
$$i_k\leftarrow i_k+1 \ \ \ \ \ak = \delta^{i_k}\alpha_{k-1} $$
and go to Step 3.
\item[5]  {Set} the new iterate $ \xkk =  \xk_+ $
\end{AlgorithmSteps}
\end{itemize}
\textsc{End}}
\end{algorithm}
It consists in a variable metric forward--backward iteration (Step 4) combined with an extrapolation--projection step (Step 1).

The steplength $\alpha_k$ is adaptively computed via a backtracking procedure, while suitable choices for the extrapolation parameter $\beta_k$ and the scaling matrix $D_k$ will be described during the convergence analysis in sections \ref{sec:rate} and \ref{sec:conv_iter}.

Algorithm \ref{SGEM} can be considered a generalization of the Fast Iterative Soft Tresholding Algorithm (FISTA, \cite{Beck-Teboulle-2009b}), whose iteration has the form \eqref{FISTA1}--\eqref{FISTA2}. The main novelties we introduce with respect to FISTA are the possibility to employ at each iteration the variable metric induced by the matrix $D_k$ and the projection of the extrapolated point $\xk+\beta_k(\xk-\xkb)$, which allows to handle problems where $\dom(f)\supseteq Y$ does not coincide with the entire space $\R^n$. When $Y=\R^n$, FISTA is recovered by setting $D_k=I$ for all $k\geq 0$.

For convenience, we restate below our hypotheses on problem \eqref{minf}, which we assume to be fulfilled throughout this section.
\begin{itemize}
\item[(A1)] $f,g:\R^n\to \R\cup\{\infty\}$ are proper, convex and lower semicontinuous;
\item[(A2)] $f$ is  differentiable with $L$-Lipschitz continuous gradient on $Y\subseteq \dom(f)$, $Y$ is closed and convex and $ \dom(f) \supseteq Y\supseteq \dom(g) $;
\item[(A3)]$g$ is bounded from below over its domain;
\item[(A4)]problem \eqref{minf} admits at least a solution.
\end{itemize}
Moreover, we will indicate hereafter by $\{\xk\}$ the sequence generated by Algorithm \ref{SGEM}, while $\xs$ will denote any of the solutions o \eqref{minf}.

First, we observe that, {thanks to assumption (A2), Algorithm \ref{SGEM} is well defined, i.e. the backtracking loop between Steps {3 and 4} terminates in a finite number of steps. Indeed, from Lemma \ref{cor:1}, observing that the sequence $\{\ak\}$ is {non--increasing} and that the reducing factor is $\delta<1$, we obtain the following inequalities
\begin{equation}\label{ine_ak}
0<\frac{\delta\eta}{L}\leq \alpha_k\leq \alpha_{k-1}\leq \alpha_0
\end{equation}
In particular, the backtracking condition implies that the new iterate $\xkk$ satisfies
\begin{equation}\label{ine1}
f(\xkk) \leq \ell(\xkk;\yk) + \frac 1{2\ak}\|\yk-\xkk\|^2_\Dk
\end{equation}

In the following sections we will show that Algorithm \ref{SGEM} with a proper parameters setting has a ${\mathcal O}(1/k^2)$ convergence rate with respect to the objective function values, i.e.
\begin{equation}\nonumber F(\xk) -F(\xs) = {\mathcal O}\left(\frac 1 {k^2}\right)\end{equation}
which is the same as FISTA, and, moreover, the iterate sequence $\{\xk\}$ converges to a minimizer of \eqref{minf}.
\subsection{Convergence rate analysis}\label{sec:rate}
In the rest of the paper we will assume that the extrapolation parameter $\beta_k$ has the form
\begin{equation}\label{betak}
\beta_k =   \frac{\thk(1-\thkb)}{\thkb}
\end{equation}
where $\{\thk\}\subset (0,1]$ is a given sequence of parameters. Moreover, we will adopt the following notation
\begin{eqnarray}
v_k &=& F(\xk)-F(\x^*) \nonumber \\
\zk& =& \xk+ \frac{1-\thkb}{\thkb}(\xk-\xkb) = \xkb +\frac{1}{\thkb} (\xk-\xkb) \label{zeta}\\
\uk &=& \zk-\x^* \nonumber\\
\tk &=& \frac{1}{\theta_k} \nonumber
\end{eqnarray}
Before giving the main result, we need to prove some technical lemmas. The first of them establishes a key inequality which is crucial for the subsequent analysis and it is analogous to Lemma 4.1 in \cite{Beck-Teboulle-2009b}.
\begin{Lemma}\label{lemma:3}
Let $\{D_k\}\subset\D_\eta$ be a sequence of scaling matrices and assume that
$\{\thk\}$ satisfies
\begin{equation}\label{thetak1}
\frac{1-\thk}{\thk^2}\leq \frac{1}{\theta_{k-1}^2} \quad \quad 0 < \thk\leq 1
\end{equation}
Then, we have
\begin{equation}\label{ricorsiva}
2\akk \tk^2\vkk + \|\ukk\|^2_\Dk \leq 2\ak t_{k-1}^2\vk+\|\uk\|^2_\Dk
\end{equation}
\end{Lemma}
\proof Let us define the point $\y^* = (1-\thk)\xk+\thk\x^*$.  {We have $\y^*\in \dom(g)$}.  {From \eqref{lemma:2d} in Lemma \ref{lemma:2c} with $y=\xkk$, $x=\yk$ and $z=\y^*$, we have
\begin{eqnarray*}
&& F(\xkk) + \frac 1{2\ak}\|\y^*-\xkk\|^2_\Dk \leq F(\y^*) +\frac 1 {2\ak}\|\y^*-\yk\|^2_\Dk\\
&&\quad \leq (1-\thk) F(\xk) +\thk F(\x^*) + \frac 1{2\ak}\|\y^*-\yk\|^2_\Dk \\
&&\quad \leq (1-\thk) F(\xk) +\thk F(\x^*) + \frac 1{2\ak}\|\y^*-(\xk+\beta_k(\xk-\xkb))\|^2_\Dk \\
\end{eqnarray*}
where the first inequality follows from the definition of $\y^*$ and the convexity of $F$ and the second one from {\eqref{nonexpansive}}. {From \eqref{betak} and the definition of $\y^*$, we obtain}
\begin{eqnarray*}
 && F(\xkk) + \frac 1{2\ak}\|(1-\thk)\xk+\thk\x^*-\xkk\|^2_\Dk  \leq (1-\thk) F(\xk) +\thk F(\x^*) + \\
& &\quad \quad + \frac 1{2\ak}\left\|(1-\thk)\xk+\thk\x^*-\left(\xk+\frac{\thk(1-\thkb)}{\thkb}(\xk-\xkb)\right)\right\|^2_\Dk
\end{eqnarray*}
{Rearranging terms, we have
\begin{eqnarray*}
&& F(\xkk) + \frac{\thk^2}{2\ak}\|\zkk-\x^*\|_\Dk^2 \leq  (1-\thk) F(\xk) +\thk F(\x^*) + \frac {\thk^2}{2\ak}\|\zk-\x^*\|_\Dk^2
\end{eqnarray*}}
Subtracting $F(\x^*)$ from both sides leads to
\begin{eqnarray*}
\vkk  + \frac{\thk^2}{2\ak}\|\zkk-\x^*\|_\Dk^2 
&\leq&  (1-\thk) \vk  + \frac {\thk^2}{2\ak}\|\zk-\x^*\|_\Dk^2
\end{eqnarray*}
Multiplying both sides by $2\ak/\thk^2$ and rearranging terms gives
\begin{equation}
2\frac{\ak}{\thk^2} \vkk + \|\zkk-\x^*\|_\Dk^2 \leq 2\ak \frac{1-\thk}{\thk^2}\vk + \|\zk-\x^*\|_\Dk^2
\label{ineini-0}
\end{equation}
Finally, observing that $\akk\leq \ak$,  {we obtain
\begin{equation}
2\frac{\akk}{\thk^2} \vkk + \|\zkk-\x^*\|_\Dk^2 \leq 2\ak \frac{1-\thk}{\thk^2}\vk + \|\zk-\x^*\|_\Dk^2\label{ineini}
\end{equation}}
In view of \eqref{thetak1}, we obtain \eqref{ricorsiva}.
\endproof
An example of sequence $\{\theta_k\}$ and corresponding $\{\beta_k\}$ satisfying \eqref{betak}-\eqref{thetak1} is the following one
\begin{equation}\label{thetak2}
\thk = \left\{\begin{array}{ll}
1 & k=0\\
\frac{a}{k+a} & k\geq 1
\end{array}\right.\ \ \quad
\beta_k = \left\{\begin{array}{ll}
0&k=0\\
\frac{k-1}{k+a} & k\geq 1
\end{array}\right.
\end{equation}
with $a\geq 2$. Indeed, since $\theta_k=\frac{1}{t_k}$, {condition} \eqref{thetak1} writes also as
\begin{equation}\nonumber
t_{k-1}^2+t_k-t_k^2\geq 0
\end{equation}
Since \eqref{thetak2} implies $t_k=\frac{k+a}{a}$, for all $k\geq 0$, $a\geq 2$ we have
\begin{eqnarray*}
t_{k-1}^2+t_k-t_k^2&=&\frac{(k-1+a)^2}{a^2}+ \frac{(k+a)}{a}-\frac{(k+a)^2}{a^2}\\
&=& \frac{(k+a)(a-2)+1}{a^2}\geq 0
\end{eqnarray*}
The choice $a=2$ has been proposed in \cite{Beck-Teboulle-2009b} for computing FISTA's extrapolation parameters, while the more general case $a\geq 2$ is considered in \cite{Chambolle-Dossal-2014}.

Our aim is now to show that the sequence $\{\vk\}$ is bounded. 
To this end, we recall the following lemma on summable nonnegative sequences.
\begin{Lemma}\label{lemma:4-1}\cite{Polyak-1987}
Let $\{a_k\}$, $\{\zeta_k\}$ and $\{\xi_k\}$ be nonnegative sequences of real numbers such that $a_{k+1} \leq (1+\zeta_k)a_k + \xi_k$ and $\sum_{k=0}^\infty \zeta_k < \infty$, $\sum_{k=0}^\infty \xi_k < \infty$. Then, $\{a_k\}$ converges.
\end{Lemma}
In the next lemma, we introduce {a crucial assumption on the sequence of matrices $\{D_k\}$ (see also \cite{Combettes-Vu-2013,Combettes-Vu-2014}).
\begin{Lemma}\label{lemma:6}
Let $\{\theta_k\}$ satisfy \eqref{thetak1} and define {$a_k=2\ak t_{k-1}^2\vk+\|\uk\|_{D_{k}}^2$}. Assume that the sequence of matrices $\{D_k\}\subset \D_\eta$,  satisfies
\begin{equation}\label{Dkipo}
D_{k+1}\preceq (1+\eta_{k})D_k\quad \forall k\geq 0\quad  \mbox{with }\eta_k\in \R, \eta_k\geq 0 \mbox{ such that }  \sum_{k=0}^\infty \eta_k < \infty
\end{equation}
Then, $\{a_k\}$ is a convergent sequence.
\end{Lemma}
\proof Setting $\sk$ as
\begin{equation}\label{sk}\sk=2\ak t_{k-1}^2\vk\end{equation}
in view of \eqref{Dkipo}, we obtain
\begin{eqnarray*}
a_{k+1}= \skk + \|\ukk\|^2_{D_{k+1}} &\leq & \skk + (1+\eta_k) \|\ukk\|^2_{D_{k}}\\
&\leq & (1+\eta_k) (\skk + \|\ukk\|^2_{D_{k}}) \\
&\leq& (1+\eta_k) (\sk +  \|\uk\|^2_{D_{k}})\\
&=& (1+\eta_{k}) a_k
\end{eqnarray*}
where the second inequality follows from the fact that $\eta_{k}\geq 0$ for any $k\geq 0$ and the third one from inequality \eqref{ricorsiva}.
Thus, Lemma \ref{lemma:4-1} implies that $\{a_k\}$ converges. \endproof
{We are now ready to give the main result of this section, establishing the convergence rate of Algorithm \ref{SGEM}. }
\begin{Thm}\label{thm:1}
Let $\{D_k\}\subset \D_\eta$ be a sequence of matrices satisfying \eqref{Dkipo} and assume that $\{\theta_k\}$, $\{\beta_k\}$ are chosen as in \eqref{thetak2} with $a\geq 2$. Then, there exists a constant $C$ such that
\begin{equation}\label{rate}F(\xk)-F(\x^*)\leq \frac{C}{(k-1+a)^2} \end{equation}
for all $k\geq 1$.
\end{Thm}
\proof
Lemma \ref{lemma:6} guarantees in particular that there exists a constant $K>0$ such that
$ a_k= 2 \ak t_{k-1}^2\vk + \|\uk\|^2_{D_{k}}\leq K $.
Since $2\ak t_{k-1}^2\vk\leq a_k$, we also have $ 2 \ak t_{k-1}^2\vk\leq K$. Formula \eqref{thetak2} implies that
{$
t_{k-1}^2=\frac{(k-1+a)^2}{a^2}$.
Thus, recalling the definition of $\vk$ and the lower bound in \eqref{ine_ak} for the parameter $\ak$, we can write
$$\vk=F(\xk)-F(\x^*)\leq \frac{a^2 L K}{2\eta\delta(k-1+a)^2} $$
{obtaining \eqref{rate} with $C=\frac{a^2LK}{2\eta\delta}$.}
\endproof
Theorem \ref{thm:1} is a generalization of Theorem 4.4 in \cite{Beck-Teboulle-2009b}, which is recovered when $D_k=I$, $a=2$, $Y=\R^n$. An analogous result for the case $D_k=I$, $a\geq 2$ can be found in \cite{Chambolle-Dossal-2014}.

It can be observed that the optimal value of the constant $C$ in \eqref{rate} is obtained with $a=2$. However, as pointed out in \cite{Chambolle-Dossal-2014} and as we will see in the following section, selecting $a > 2$ allows to prove the convergence of the iterates $\{\xk\}$ to a solution of \eqref{minf}.
\paragraph{Remark}When {the sequence of matrices $\{D_k\}$ satisfies the assumption \eqref{Dkipo} and, in addition, $\mbox{sup}_{k\in \N}\|D_k\|=\gamma<\infty$, then there exists a matrix $\mathfrak{D}\in \D_{\eta}$ such that $D_k\rightarrow \mathfrak{D}$ pointwise \cite[Lemma 2.3]{Combettes-Vu-2013}. A similar result holds also for a sequence of matrices $\{D_k\}\subset\D_\eta$, such that
$\mbox{sup}_{k\in \N}\|D_k\|=\gamma<\infty$, {satisfying the following condition} 
\begin{equation}\label{Dk1}
D_{k}\preceq (1+\nu_{k})D_{k+1}\quad \forall k\geq 0\quad  \mbox{with }\nu_k\in \R,\ \nu_k\geq 0 \mbox{ such that }  \sum_{k=0}^\infty \nu_k < \infty
\end{equation}
A sufficient condition ensuring both \eqref{Dkipo} and \eqref{Dk1} is the following one
\begin{equation}\label{choice}
 \frac{1}{\gamma_k} \leq \|D_k\|\leq \gamma_k \quad \gamma_k^2=1+\zeta_k \quad \mbox{where } \ \ \ \zeta_k \geq 0\ \ \ \mbox{ and } \ \ \ \sum_{k=0}^\infty \zeta_k<\infty
 \end{equation}
 with $\gamma_k<\gamma$, $\gamma>1$. Indeed, in this case $\eta=\frac{1}{\gamma}$ and for any $x\in \R^n$ we have
 \begin{eqnarray*}
  x^T D_{k+1}x&\leq & \frac{\gamma_{k+1}\gamma_k}{\gamma_k}\|x\|^2 \leq  \gamma_{k+1}\gamma_k x^T D_k x\\
  x^T D_{k}x&\leq & \frac{\gamma_{k}\gamma_{k+1}}{\gamma_{k+1}}\|x\|^2 \leq  \gamma_k\gamma_{k+1} x^T D_{k+1} x\\
  \end{eqnarray*}
 Let us define $\eta_k=\nu_k=\gamma_{k+1}\gamma_{k} -1 = \sqrt{(1+\zeta_{k+1})(1+\zeta_k)}-1$ and observe that the series $\sum \eta_k$ and $\sum \zeta_k$ have the same behavior, since the known limit $\lim_{z\rightarrow 0} (\sqrt{1+z}-1)/z=1/2$. Therefore, for any $x\in \R^n$, we can write
\begin{eqnarray*}
&& x^T D_{k+1}x\leq  (1+\eta_k) x^T D_k x \\
&& x^T D_{k}x\leq  (1+\nu_k) x^T D_{k+1} x
\end{eqnarray*}
with $\eta_k=\nu_k$ for any $k\geq 0$ and $\sum \eta_k <\infty$; {thus} a sequence of matrices chosen according to the rule \eqref{choice} satisfies both the assumptions \eqref{Dkipo} and \eqref{Dk1}.

For the sequence $\{D_k\}$ satisfying \eqref{choice} and $\{\theta_k\}$ as in \eqref{thetak2} with $a=2$, the convergence rate estimate established in Theorem \ref{thm:1} becomes
$$F(\xk)-F(\x^*)\leq \frac{2L\gamma K}{\delta (k+1)^2} $$
which is very similar to that in \cite{Beck-Teboulle-2009b}.
\subsection{Convergence of the iterates to a minimizer}\label{sec:conv_iter}
In this section, borrowing the ideas in \cite{Chambolle-Dossal-2014}, we prove that the iterates $\{\xk\}$ converge to a minimizer of $F$, when the parameters sequences are chosen as in \eqref{thetak2} with $a>2$ and the scaling matrices sequence satisfies some additional assumption. Before giving the main result, we need to prove some technical lemmas.

Under the same assumptions {of} Theorem \ref{thm:1}, we first prove the boundedness of the sequence $\{\|\uk\|\}$. To this end, we first recall an useful lemma on summable nonnegative sequences.
\begin{Lemma}\label{lemma:4}\cite{Bonettini-Prato-2015a}
Let $\{\gamma_k\}$ be a sequence of positive numbers such that $\gamma_k^2 = 1+\eta_k$, $\eta_k \geq 0$, where $\sum_{k=0}^\infty \eta_k<\infty$. Let $\tau_k = \prod_{j=0}^k\gamma_j^2$ for any $k\geq 0$. Then the sequence $\{\tau_k \}$ is  bounded.
\end{Lemma}
\begin{Lemma}\label{lemma:5}
Let $\{\theta_k\}$, $\{\beta_k\}$ be defined such that {\eqref{betak} and \eqref{thetak1} hold} and assume that the sequence $\{D_k\}\subset \D_{\eta}$
satisfies {{\eqref{Dkipo}}}.
Then, the sequence $\{\|\uk\|\}$ is bounded, i.e. $\|\uk\|\leq U$ for $k\geq 0$.
\end{Lemma}
Recalling definition \eqref{sk}, we obtain
\begin{eqnarray*}
\|\ukk\|_{D_{k+1}}^2&\leq& (1+\eta_k) \|\ukk\|_{D_{k}}^2 \\
                    &\leq& (1+\eta_k)  (\sk-\skk +\|\uk\|_{D_k}^2)\\
          &\leq& (1+\eta_{k}) (\sk + \|\uk\|_{D_{k}}^2)\\
          &\leq& (1+\eta_{k})(\sk+ (1+\eta_{k-1})\|\uk\|_{D_{k-1}}^2) \\
          &\leq& (1+\eta_{k})  (1+\eta_{k-1}) \left(\sk + \left( s_{k-1}-\sk+\|\u^{(k-1)}\|_{D_{k-1}}^2\right)\right)\\
          &\leq&  (1+\eta_{k}) (1+\eta_{k-1}) (s_{k-1}+ (1+\eta_{k-2})\|\u^{(k-1)}\|_{D_{k-2}}^2)\\
          &\leq&   (1+\eta_k)   (1+\eta_{k-1}) (1+\eta_{k-2}) (s_{k-2}+\|\u^{(k-2)}\|_{D_{k-2}}^2) \\
          &\leq&  \tau_k (s_0+\|u^{(0)}\|_{D_0}^2)
\end{eqnarray*}
where $\tau_k$ is defined as in Lemma \ref{lemma:4} and we repeatedly applied the inequalities \eqref{Dkipo} and \eqref{ricorsiva} in the following ones, together with the fact that $\eta_i\geq 0$ for $i\geq 0$. Since $\tau_k$ is bounded, $\|\ukk\|_{D_{k+1}}^2$ is bounded. Furthermore, from \eqref{ine_norm} we have the result. \endproof
The following lemma generalizes the results in \cite[Theorem 2]{Chambolle-Dossal-2014}.
\begin{Lemma}\label{Lemmanuova1}
Let $\{D_k\}\subset \D_\eta$ be a sequence of matrices satisfying \eqref{Dkipo}, with $\sup_{k\in \N}\|D_k\|=\gamma<\infty$.
Assume that $\{\thk\}$, $\{\beta_k\}$ are chosen as in \eqref{thetak2}, with $a> 2$. Then the sequence $\{k v_k\}$ is summable.
\end{Lemma}
\proof First we recall that with these settings for $\theta_k$ and $\beta_k$, \eqref{thetak1} and \eqref{betak} are satisfied. In view of $\tk^2 = \frac{1}{\thk^2}$,
we can write the inequality {\eqref{ineini}} as follows:
\begin{eqnarray*}
\alpha_{k+1}\tk^2 \vkk-\alpha_k (\tk^2-\tk)\vk&\leq & \frac{\|\uk\|_{\Dk}^2}{2}- \frac{\|\ukk\|_\Dk^2}{2}
\end{eqnarray*}
Summing up from $k=0,...,K$, since $t_0=1$,  we have
\begin{align*}
\alpha_{K+1}t_K^2 v_{K+1} &+\sum_{k=1}^K \alpha_k (\tkm^2-\tk^2+\tk)\vk \leq \\
&\leq \frac 1 2\sum_{k=1}^K {\|\uk\|_{\Dk}^2} - {\|\uk\|_{D_{k-1}}^2}+\frac{\|\u^{(0)}\|_{D_0}^2}{2 }- \frac{\|\u^{(K+1)}\|_{D_K}^2}{2 }\\
&\leq \frac 1 2 \sum_{k=1}^K (1+\eta_{k-1}) {\|\uk\|_{D_{k-1}}^2} - {\|\uk\|_{D_{k-1}}^2}
+\frac{\|\u^{(0)}\|_{D_0}^2}{2}\\
&= \frac{1}{2}\sum_{k=1}^K  \eta_{k-1} \|\uk\|_{D_{k-1}}^2+\frac{\|\u^{(0)}\|_{D_0}^2}{2}\\
&\leq \frac{\gamma {U^2}}{2}\sum_{k=0}^{K-1} \eta_k  +\frac{\|\u^{(0)}\|_{D_0}^2}{2}
\end{align*}
where the second inequality follows from {\eqref{Dkipo}} and the last one from the boundedness of $\|D_k\|$ and Lemma \ref{lemma:5}.
Furthermore, using \eqref{ine_ak}, we obtain
$$\sum_{k=1}^K (\tkm^2-\tk^2+\tk)\vk \leq \frac{L\gamma {U^2}}{2\delta\eta}\sum_{k=0}^{K-1} \eta_k  +\frac{\|\u^{(0)}\|_{D_0}^2}{2}$$
Then, since $\eta_k$ is a summable sequence, in view of \eqref{thetak1}, $\{(\tkm^2-\tk^2+\tk)\vk\}$ is {nonnegative} and summable. Finally, observing that for $a>2$, we have
$$0\leq\tkm^2-\tk^2+\tk= \frac{k(a-2)+(a-1)^2}{a^2} $$
we can conclude that also $\{kv_k\}$ is summable. \endproof
The following lemma is a consequence of the previous one. It requires an additional condition on the scaling matrix sequence $\{D_k\}$. Indeed we will assume that the sequence $\{\eta_k\}$ in \eqref{Dkipo} is given by
\begin{equation}\label{ketak}
\{\eta_k\}=\mathcal{O}\left(\frac{1}{k^p}\right) \quad \mbox{with }p>2
\end{equation}
This assumption guarantees also that $\{k\eta_k\}$ is summable.
When $\{D_k\}$ is chosen according to \eqref{choice}, the condition \eqref{ketak} is satisfied when $\zeta_k=\frac{b}{k^{p}}$  for any positive scalar $b$  and $p>2$.
\begin{Lemma}\label{Lemmanuova2}
{Let the assumption of Lemma \ref{Lemmanuova1} be fullfilled with
$\{\eta_k\}=\mathcal{O}(\frac{1}{k^p})$, $p>2$ in \eqref{Dkipo}.}
{Then, setting $\delta_k= {\|\xk-\xkm\|_{D_{k}}^2}/{2}$,}  the sequence $\{k \delta_k\}$ is summable. In addition, there exists $D>0$ such that for all $k\geq 1$, $\delta_k\leq \frac{D}{k^2}$.
\end{Lemma}
\proof From \eqref{lemma:2d} with $\x=\yk$, $\y=\xkk$ and $\z=\xk$, it follows that
\begin{equation}
F(\xkk) +\frac{\|\xk-\xkk\|_\Dk^2}{2\ak}\leq F(\xk) +\frac{\|\xk-\yk\|_\Dk^2}{2\ak}
\label{1}\end{equation}
From Lemma \ref{lemma:nonexpansive}, since $\xk\in \dom(g)\subseteq Y$, we have
$$\|\xk-\yk\|_\Dk^2\leq \beta_k^2 \|\xk-\xkm\|_\Dk^2$$
Then, subtracting $F(\xs)$ from both sides of \eqref{1}, we can write
\begin{equation}
\vkk +\frac{\|\xk-\xkk\|_\Dk^2}{2\ak}\leq \vk +\beta_k^2 \frac{\|\xk-\xkm\|_\Dk^2}{2\ak}
\label{2}
\end{equation}
From {\eqref{Dkipo}}  we have
$$\delta_{k+1}=\frac{1}{2}\|\xkk-\xk\|_{D_{k+1}}^2\leq (1+\eta_k)\frac{1}{2}\|\xkk-\xk\|_{D_{k}}^2$$
Then, since $\eta_k\geq 0$, from \eqref{2}, it follows that
\begin{eqnarray} \label{3}
\delta_{k+1}&\leq&
 (1+\eta_k) (\alpha_k (\vk -\vkk)+\beta_k^2 \delta_k )
%
\end{eqnarray}
{Since  $\theta_k=\frac{a}{k+a}$ and $\beta_k=\frac{k-1}{k+a}$, \eqref{3} writes also as
\begin{eqnarray}
(k+a)^2 \delta_{k+1} - (1+\eta_k) (k-1)^2 \delta_k &\leq&  (1+\eta_k)\alpha_k (k+a)^2(\vk -\vkk)\nonumber \\
&\leq & \alpha_0 (1+\eta_k) (k+a)^2 (\vk -\vkk) \label{4}
\end{eqnarray}
where the second inequality follows from  \eqref{ine_ak}.

Since $\{\eta_k\}=\mathcal{O}(\frac{1}{k^p})$ with $p>2$,  $\lim_{k\rightarrow \infty} \eta_k k^2=0$; therefore,  given a scalar $\epsilon>0$ such that
$a^2-2a\geq \epsilon$, there exists an index $\ell$ such that for any {$k> \ell$} we can write
\begin{equation}\label{6}
\eta_k (k-1)^2<\eta_k k^2< \epsilon\leq a^2-2a
\end{equation}
Summing up the inequality \eqref{4} for $k=\ell,...,K$ yields
\begin{eqnarray}
&& (K+a)^2 \delta_{K+1}+ \sum_{k=\ell+1}^K ((k-1+a)^2-(1+\eta_k)(k-1)^2)\delta_k  \leq (1+\eta_{\ell})(\ell-1)^2 \delta_{\ell} + \label{5_0}\\
&& \quad +\alpha_0 ( (\ell+a)^2 (1+\eta_{\ell})v_{\ell}- (K+a)^2(1+\eta_K)) v_{K+1}+  \nonumber\\
 && \quad +\alpha_0 \sum_{k=\ell+1}^K ((k+a)^2(1+\eta_k) -(k-1+a)^2(1+\eta_{k-1})) v_k \nonumber
\end{eqnarray}
For all $k$ we have
\begin{eqnarray*}(k+a)^2(1+\eta_k) -(k-1+a)^2(1+\eta_{k-1}) &=& \eta_k(k+a)^2+2(k+a)-1-(k-1+a)^2\eta_{k-1} \\
&\leq& \eta_k(k+a)^2+2(k+a)
\end{eqnarray*}
and, in view of \eqref{6}, for $k>\ell$ we also have
 $$(k-1+a)^2-(1+\eta_k)(k-1)^2= a^2 +2ka -2a -\eta_k (k-1)^2 > 2ka >0$$
 Then, ignoring negative terms on the right hand side, we obtain
 \begin{eqnarray}
&&  (K+a)^2 \delta_{K+1}+ \sum_{k=\ell+1}^K 2k a \delta_k  \leq   (1+\eta_{\ell})(\ell-1)^2 \delta_{\ell}+ \label{5}\\
&&\quad +\alpha_0 \left( (\ell+a)^2(1+\eta_{\ell}) v_{\ell}+ \sum_{k=\ell+1}^K 2(k+a) v_k +\eta_k (k+a)^2 v_k\right) \nonumber
\end{eqnarray}
Since $\eta_k (k+a)^2$ is bounded, by  Lemma \ref{Lemmanuova1}, the right hand side of the previous inequality is uniformly bounded independently on $K$. This  ensures that $k\delta_k$ is summable. Furthermore $K^2 \delta_{K+1}$ is globally bounded. \endproof
We are now able to prove the first, weak, convergence result about the sequence generated by Algorithm \ref{SGEM}, as stated in the following Corollary.
\begin{Cor}\label{cor:2}
Let the assumptions of Lemma \ref{Lemmanuova2} be satisfied. Then, $\{\xk\}$ is bounded and any of its limit point is a solution of problem \eqref{minf}.
\end{Cor}
\proof A direct consequence Lemma \ref{Lemmanuova2} is that the sequence $\{k(\xk-\xkm)\}$ is bounded. From Lemma \ref{lemma:5}, it follows that the sequence $\{\zk\}$ defined in \eqref{zeta} is also bounded. These two facts imply that the sequence $\{\xk\}$ is bounded.

Assume that $\tilde{x}$ is a limit point of $\{\xk\}$, i.e. there exists a subsequence $\{\xk\}_{k\in \mathcal{K}}$ of $\{\xk\}$ such that $\xk\rightarrow \tilde{x}$, $k\in \mathcal{K}$ as $k\rightarrow\infty$. This element $\tilde{x}$ of $\dom(g)$ is also a limit point of $\{\yk\}$. 
Indeed, from Lemma \ref{lemma:nonexpansive}, the definition of $\yk$ and the boundedness of $\beta_k$ and $\{D_k\}$, we have for any $k\geq 1$ and, in particular for $k\in \mathcal{K}$:
$$\|\yk-\xk\|^2\leq \frac{1}{\eta} \|\yk-\xk\|^2_{D_k}\leq {\frac{2}{\eta}}\beta_k^2 \delta_k\leq {\frac{2}{\eta}} \delta_k$$
Under the assumption of Lemma \ref{Lemmanuova2}, we have $\delta_k\rightarrow 0$, as $k\rightarrow \infty$, thus $\tilde{x}$ is a limit point of $\{\yk\}$. From assumption \eqref{Dkipo} we have that $\{D_k\}$ converges pointwise to some matrix  $\mathfrak{D}\in \D_{\eta}$ \cite[Lemma 2.3]{Combettes-Vu-2013} and, since $\alpha_k\in [\delta\eta/L,\alpha_0]$ we can assume without loss of generality that $\alpha_k$ converges to some $\alpha>0$ as $k$ diverges, $k\in \mathcal{K}$. Therefore, $\tilde{x}$ is a fixed point of the operator $p_{\alpha, \mathfrak{D}}$ and, consequently, it is a minimizer of $F$.
\endproof
Before giving the main result stating that the whole sequence converges to a minimizer, we need to prove the following technical lemma, which holds when the matrices sequence $\{D_k\}$ satisfies both \eqref{Dkipo} and \eqref{Dk1}.
\begin{Lemma}\label{lemmaagg}
{Let the assumption of Lemma \ref{Lemmanuova1} be fullfilled and
$\{D_k\}$ be a sequence of matrices satisfying both the conditions \eqref{Dkipo} and \eqref{Dk1}. }
Then, denoting $\Phi_k=\frac{\|\xk-\xs\|_{D_k}^2}{2}$, for $k\geq 1$ we have
\begin{equation}\label{7}
\Phi_{k+1} -\Phi_k \leq \beta_k (\Phi_k -\Phi_{k-1}) +2\beta_k \delta_k (1+\eta_k) +  \eta_k (1+\beta_k) \Phi_k+ \beta_k \nu_{k-1}\frac{\|\xs-\xkm\|_{\Dk}^2}{2}
\end{equation}
\end{Lemma}
\proof Let $\xs$ be a solution of the problem \eqref{minf}. Using \eqref{lemma:2d} in Lemma \ref{lemma:2c} with $\y=\xkk$, $\x=\yk$, $\z=\xs$, we have
\begin{align}
F(\xkk)&+\frac{1}{2\alpha_k}\|\xs-\xkk\|_{\Dk}^2 \leq F(\xs)+\frac{1}{2\alpha_k}\|\xs-\yk\|_{\Dk}^2\nonumber\\
&\leq F(\xs)+\frac{1}{2\alpha_k}\|\xs-\xk\|_{\Dk}^2+\frac{\beta_k^2}{2\alpha_k}\|\xk-\xkm\|_{\Dk}^2 +\frac{\beta_k}{\alpha_k} (\xk-\xkm)^T \Dk (\xk-\xs)\label{20}
\end{align}
where the last inequality follows from  Lemma \ref{lemma:nonexpansive}, since $\yk=P_{Y,\Dk}(\xk+\beta_k(\xk-\xkm))$.
We observe that
\begin{eqnarray*}
&& (\xk-\xkm)^T \Dk (\xk-\xs)=  
\frac{1}{2}\|\xk-\xkm\|_\Dk^2+ \frac{1}{2}\|\xk-\xs\|_\Dk^2 - \frac{1}{2}\|\xkm-\xs\|_\Dk^2
\end{eqnarray*}
Using this equality in \eqref{20} we obtain
\begin{eqnarray*}
&&F(\xkk)+\frac{1}{2\alpha_k}\|\xs-\xkk\|_{\Dk}^2 \leq F(\xs)+\frac{1}{2\alpha_k}\|\xs-\xk\|_{\Dk}^2+\frac{\beta_k^2}{2\alpha_k}\|\xk-\xkm\|_{\Dk}^2+ \nonumber\\
&& +\frac{\beta_k}{\alpha_k} \left( \frac{1}{2}\|\xk-\xkm\|_\Dk^2+ \frac{1}{2}\|\xk-\xs\|_\Dk^2 - \frac{1}{2}\|\xkm-\xs\|_\Dk^2\right)\label{21}
\end{eqnarray*}
Then we have
\begin{eqnarray*}
&&\frac{1}{2\alpha_k}\|\xs-\xkk\|_{\Dk}^2 -\frac{1}{2\alpha_k}\|\xs-\xk\|_{\Dk}^2 \leq F(\xs)-F(\xkk) +\frac{\beta_k^2+\beta_k}{2\alpha_k}\|\xk-\xkm\|_{\Dk}^2+ \nonumber\\
&& +\frac{\beta_k}{\alpha_k} ( \frac{1}{2}\|\xk-\xs\|_\Dk^2 - \frac{1}{2}\|\xkm-\xs\|_\Dk^2)\label{21a}
\end{eqnarray*}
Since $F(\xs)-F(\xkk) \leq 0$ and $\beta_k^2+\beta_k\leq 2\beta_k$, we obtain
\begin{equation*}
 \frac{1}{2}\|\xs-\xkk\|_{\Dk}^2 - \frac{1}{2}\|\xs-\xk\|_{\Dk}^2 \leq  {\beta_k} (\frac{1}{2}\|\xs-\xk\|_\Dk^2 - \frac{1}{2} \|\xs-\xkm\|_\Dk^2) + 2 \beta_k \delta_k
\end{equation*}
Multiplying the last inequality by $(1+\eta_k)$, from \eqref{Dkipo}, we obtain 
\begin{align*}
\Phi_{k+1}&- (1+\eta_k)\Phi_k \leq  \\
& \leq  {\beta_k}(1+\eta_k) (\Phi_k - \frac{\|\xs-\xkm\|_\Dk^2}{2}) + 2 \beta_k (1+\eta_k)\delta_k\\
& = \beta_k (\Phi_k-\Phi_{k-1}) + \beta_k \eta_k \Phi_k +\beta_k (\Phi_{k-1}- (1+\eta_k)\frac{\|\xs-\xkm\|_\Dk^2}{2})+ 2 \beta_k (1+\eta_k)\delta_k
\end{align*}
Thus, in view of the assumption \eqref{Dk1}, since
$$\Phi_{k-1}=\frac{\|\xs-\xkm\|^2_{D_{k-1}}}{2}\leq (1+\nu_{k-1}) \frac{\|\xs-\xkm\|^2_{D_{k}}}{2}$$
we obtain \eqref{7}.
\endproof
Now, {as in \cite{Chambolle-Dossal-2014}}, we introduce the notation
\begin{eqnarray}
\beta_{j,k}&=&\Pi_{\ell=j}^k \beta_{\ell}=\Pi_{\ell=j}^k \frac{\ell-1}{\ell+a} \quad j\geq 1, k\geq j \nonumber\\
\beta_{j,k}&=&1 \quad j>k \label{defbjk}
\end{eqnarray}
Since $\beta_1=0$, $\beta_{1,k}=0$ for $k\geq 1$. Moreover, in  \cite{Chambolle-Dossal-2014} it is proved that,
for $a>2$, the following inequality {holds} for $j\geq2$
\begin{equation}\label{CD}
\sum_{k=j}^\infty \beta_{j,k}\leq \frac{j+5}{2}
\end{equation}
{This inequality is exploited in the proof of the following convergence theorem, whose line is very similar to that of Theorem 3 in \cite{Chambolle-Dossal-2014}. This result requires that the sequence of matrices $\{D_k\}$ satisfies both the assumptions \eqref{Dkipo} and \eqref{Dk1}, where $\{\eta_k\}$ and $\{\nu_k\}$ are $\mathcal{O}(\frac{1}{k^p})$ with $p>2$.}
\begin{Thm}\label{convergenza}
Assume that $\{\theta_k\}$ and $\{\beta_k\}$ are chosen as in \eqref{thetak2} with $a>2$ and
let $\{D_k\}\subset \D_\eta$ be a sequence of matrices satisfying \eqref{Dkipo} and \eqref{Dk1} with $\sup_{k\in \N}\|D_k\|=\gamma<\infty$,
$\{\eta_k\}=\mathcal{O}(\frac{1}{k^p})$ and $\{\nu_k\}=\mathcal{O}(\frac{1}{k^p})$ with $p>2$.  Then,  the sequence $\xk$ converges to a minimizer of $F$.
\end{Thm}
\proof  The proof follows \cite[Theorem 3]{Chambolle-Dossal-2014} and \cite{Lorentz-Pock-2014}.
We first prove that $\Phi_k= \|\xk-\xs\|_{D_k}^2/2$ converges.
From Corollary \ref{cor:2}, we have that $\{\xk\}$ is a bounded sequence. Then,
there exists a positive scalar $M$ such that $\frac{\|\xk-\xs\|^2}{2}\leq M$, for all $k\geq 0$.
From the inequality \eqref{7}, since $\sup_{k\in \N}\|D_k\|=\gamma$, we obtain
\begin{equation}\label{8}
\Phi_{k+1} -\Phi_k \leq \beta_k (\Phi_k -\Phi_{k-1}) +2\beta_k \delta_k (1+\eta_k) + \eta_k (1+\beta_k) \gamma M+ \beta_k \nu_{k-1} \gamma M
\end{equation}
Now, defining $p_k=\max(0,\Phi_k -\Phi_{k-1})$ and recalling that $\beta_k\leq 1$, we obtain
\begin{equation}\label{9}
p_{k+1} \leq \beta_k p_k +2\beta_k \delta_k (1+\eta_k) +  \beta_k \gamma M (\eta_k + \nu_{k-1}) + \gamma M \eta_k
\end{equation}
By applying \eqref{9} recursively and using \eqref{defbjk} and $\beta_1=0$, it follows that
\begin{equation}\nonumber
p_{k+1} \leq 2 \sum_{j=2}^k \beta_{j,k} \delta_j (1+\eta_j) + \gamma M \sum_{j=2}^k \beta_{j,k} (\eta_j + \nu_{j-1}) + \gamma M  \sum_{j=2}^{k} \beta_{j,k}  \eta_{j-1} + \gamma M \eta_k
\end{equation}
for all $k\geq 2$. Hence,
\begin{eqnarray*}
\sum_{k=2}^{+\infty}p_{k+1} &\leq &  2 \sum_{k=2}^{+\infty} \sum_{j=2}^k \beta_{j,k} \delta_j (1+\eta_j) + \gamma M \sum_{k=2}^{+\infty} \sum_{j=2}^k \beta_{j,k} (\eta_j + \nu_{j-1}) + \nonumber\\
&&+ \gamma M  \sum_{k=2}^{+\infty} \sum_{j=2}^{k} \beta_{j,k}  \eta_{j-1} +  \gamma M \sum_{k=2}^{+\infty}\eta_k\\
&\leq&  2 \sum_{j=2}^{+\infty}  \delta_j (1+\eta_j) \sum_{k=j}^{+\infty}  \beta_{j,k} + \gamma M \sum_{j=2}^{+\infty} (\eta_j + \nu_{j-1}) \sum_{k=j}^{+\infty}  \beta_{j,k} + \nonumber\\
&&+ \gamma M  \sum_{j=2}^{+\infty} \eta_{j-1} \sum_{k=j}^{+\infty}  \beta_{j,k}   + \gamma M \sum_{k=2}^{+\infty}  \eta_k \nonumber\\
&\leq&  2 \sum_{j=1}^{+\infty} \delta_j (1+\eta_j) \frac{j+5}{2} + \gamma M \sum_{j=1}^{+\infty} (\eta_j + \nu_{j-1})  \frac{j+5}{2}   + \nonumber\\
&&+ \gamma M  \sum_{j=1}^{+\infty} \eta_{j-1} \frac{j+5}{2}   + \gamma M \sum_{k=1}^{+\infty}  \eta_k \label{12}
\end{eqnarray*}
where the last inequality follows form \eqref{CD}.
{From the assumption on $\{\eta_j\}$ and $\{\nu_j\}$, $\{j \eta_j\}$ and  $\{j\nu_j\}$ are summable; from  Lemma \ref{Lemmanuova2},  $\{j \delta_j\}$ is also summable.
This implies that} the right side of the last inequality is finite,
therefore $\{p_k\}$ is summable. We set $q_k=\Phi_k- \sum_{i=1}^k p_i$ and since $\Phi_k\geq 0$ and $\sum_{i=1}^{\infty} p_i$ is bounded, we have that $q_k$ is bounded from below. On the other hand
$$q_{k+1}= \Phi_{k+1}-p_{k+1}-\sum_{i=1}^kp_i\leq \Phi_{k+1}-\Phi_{k+1}+\Phi_k -\sum_{i=1}^kp_i=q_k$$
Therefore $\{q_k\}$ is a non-increasing sequence and  it is convergent. This implies that $\Phi_k=s_k+\sum_{i=1}^k p_i$ is convergent.

Assume now that $\tilde \x\in \dom(g)$ is a limit point of $\{\xk\}$, i.e. there exists a subsequence $\{\x^{(k_i)}\}$ of $\{\xk\}$ such that $\lim_{i\to \infty}\x^{(k_i)}=\tilde \x$. By Corollary \ref{cor:2}, $\tilde \x$ is a minimizer of $F$. Thus, the first part of the proof applies also to $\tilde\x$ and we can conclude that $\{\|\xk-\tilde\x\|^2_{D_k}\}$ converges. As a consequence of this we have
 $$\lim_{k\to\infty}\|\xk-\tilde\x\|^2_{D_k}=\lim_{i\to\infty}\|\x^{(k_i)}-\tilde\x\|^2_{D_{k_i}}=0$$
Since $\eta \|\xk-\xs\|^2 \leq \|\xk-\xs\|^2_{D_k}$, the last equality implies that the whole sequence converges to the minimizer $\tilde\x$.
 \endproof
As a consequence of the analysis performed in this section, we can conclude that when the sequence of matrices $\{D_k\}$ is chosen so that the condition \eqref{choice} holds with $\{\zeta_k\}=\mathcal{O}\left(\frac{1}{k^p}\right)$ with $p>2$, and the extrapolation parameters $\theta_k$ and $\beta_k$ are defined as in \eqref{thetak2}, Algorithm \ref{SGEM} generates a sequence $\{\xk\}$ convergent to a minimizer of $F$ with a $\mathcal{O}\left(\frac{1}{k^2}\right)$ convergence rate for the objective function values.
\section{Numerical experiments}\label{sec:num}
In this section we present the results of several numerical experiments which aim at evaluating the {effectiveness}
of the proposed scaled forward-backward method with extrapolation (SFBEM) by comparison with other state-of-the-art algorithms. The numerical experiments concern three different optimization problems which can be formalized as in \eqref{minf} and arise from some relevant real-life applications.
\subsection{Image deblurring with Poisson noise}
We consider the inverse problem of recovering an unknown image $x_{\rm true}$ from a given data corrupted by noise. Bayesian approaches suggest to address this problem by minimizing a functional which can be expressed as the sum of a discrepancy function, typically depending on the noise type affecting the data, and a regularization term adding a priori information and possible constraints. In particular, in the case of Poisson noise, the discrepancy function measuring the distance from the data $b\in \mathbb{R}^n$ is the {generalized} Kullback-Leibler (KL) divergence of the form
\begin{equation}
\label{KL_divergence}
 \KL (x)  = \sum_{i=1}^{n} \left\{ b_i \ln\frac{b_i}{(Ax + bg)_i} + (Ax + bg)_i - b_i\right\}
\end{equation}
where $A\in\mathbb{R}^{n\times n}$ is a linear operator modeling the distortion due to the image acquisition system and $bg \in \mathbb{R}^n$ is a known {positive} background radiation constant. A typical assumption for the matrix $\ve{A}$ is that it has nonnegative elements and each row and column has at least one positive entry. We refer the interested reader to \cite{Bertero-etal-2009} for a detailed survey on the image deblurring problem in presence of Poisson noise and the properties of the KL function \eqref{KL_divergence}.

As for the regularization term, we consider a smooth discrete version of the total variation, also known in the literature as {\em hypersurface potential} (HS) \cite{AcarVog04, Charbonnier-etal-2011, Zanella-etal-2009}, {that, for a square $m\times m$ image with $m^2=n$,  is} defined as
\begin{equation}
\nonumber
\HS(x) = {\sum_{i,j=1}^m} \sqrt{((\mathcal{D} \ve{x})_{i,j})_1^2+((\mathcal{D} \ve{x})_{i,j})_2^2+\delta^2},
\end{equation}
where {$\mathcal{D}: \mathbb{R}^{m^2}\longrightarrow\mathbb{R}^{2m^2}$} is the discrete gradient operator with periodic boundary conditions
\begin{equation}\label{discr_grad}
(\mathcal{D} \ve{x})_{i,j} =
\left(
\begin{array}{c}
 ((\mathcal{D}\ve{x})_{i,j})_1\\
 ((\mathcal{D}\ve{x})_{i,j})_2
\end{array}
\right)=
\left(
\begin{array}{c}
 x_{i+1,j}-x_{i,j}\\
 x_{i,j+1}-x_{i,j}\\
\end{array}
\right), \ \ \ x_{n+1,j} = x_{1,j}, \ \ \ x_{i,n+1} = x_{i,1}.
\end{equation}
In conclusion, a way to recover the true image from the corrupted data is to find a solution of the following optimization problem
\begin{equation}
\label{min_KL+HS}
 \min_{\ve{x}\in\mathbb{R}^n} F(\ve{x}) \equiv \KL(\ve{x}) + \rho\HS(\ve{x}) + \iota_{\ve{x}\geq \ve{0}}(\ve{x}),
\end{equation}
where $\rho$ is a positive parameter balancing the role of the regularization term and $\iota_{\ve{x}\geq \ve{0}}$ denotes the indicator function of the {nonnegative} orthant; {indeed, the unknown (the pixels of the image) have to be  {nonnegative}}.

Problem \eqref{min_KL+HS} can be naturally cast in the form \eqref{minf} by setting $f(x) = \KL(x)+\rho\HS(x)$ and $g(\x) = \iota_{\ve{x}\geq \ve{0}}(\ve{x})$.

In this case $\dom(f) = \{\x\in\R^n: A\x + bg > 0\}$ and $\nabla f$ is Lipschitz continuous on $Y=\{\x\in\R^n:\x\geq 0\} =\dom(g)\subseteq\dom(f)$. However, only an estimation from above of the Lipschitz constants of both $\nabla\hspace{-0.1cm}\KL$ and $\nabla\hspace{-0.1cm}\HS$ is known ({see} \cite{Harmany12} and \cite{Jensen12}, respectively).

The numerical tests have been performed by solving the optimization problem \eqref{min_KL+HS} on two different datasets. The original images are the $128\times 128$ micro \cite{Willett-Nowak-2003} and the $256\times256$ Cameraman, both used in several papers and reported in the first row of Figure \ref{Immagini_KL+HS}. The values of the micro original image are in the range $[0,69]$, while the values of the Cameraman lay in the interval $[0,1000]$. The corrupted data (Figure \ref{Immagini_KL+HS}, second row) have been generated by convolving the objects with a {suitable} point spread function (the psf proposed in \cite{Willett-Nowak-2003} for the micro image and a Gaussian psf, with standard deviation equal to 1.3, for the Cameraman one), adding a constant background equal to 1 and perturbing the result of the convolution with Poisson noise (simulated through the \verb+imnoise+ function of the Matlab Image Processing Toolbox). For both test problems we assume periodic boundary conditions, thus $A$ is block-circulant with circulant blocks and the matrix vector products involving $A$ can be performed via the Fast Fourier Transform.

We chose the regularization parameter $\rho$ equal to $0.09$ for the first test problem and $0.045$ for the second one and the parameter $\delta$ for the HS functional equal to $0.05$ for both datasets.
\begin{figure}[htb!]
\begin{center}
 \begin{tabular}{cc}
\includegraphics[width=.25\textwidth]{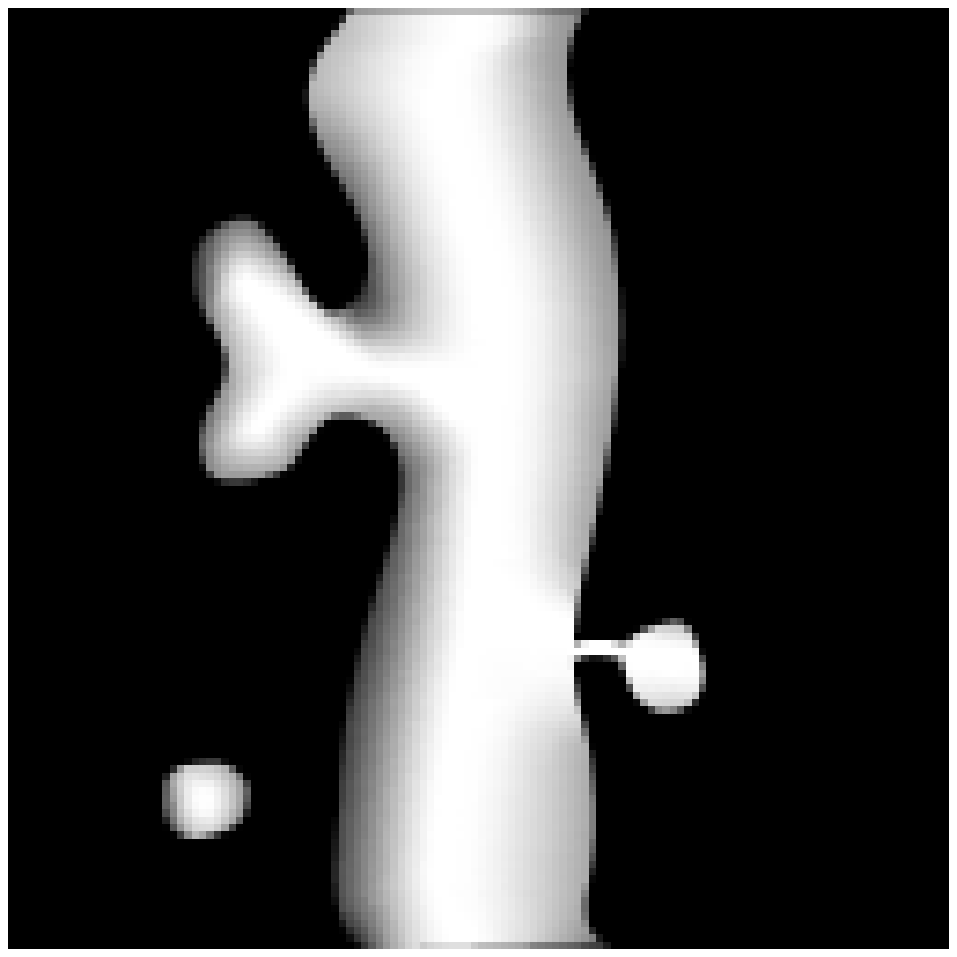}&
\includegraphics[width=.25\textwidth]{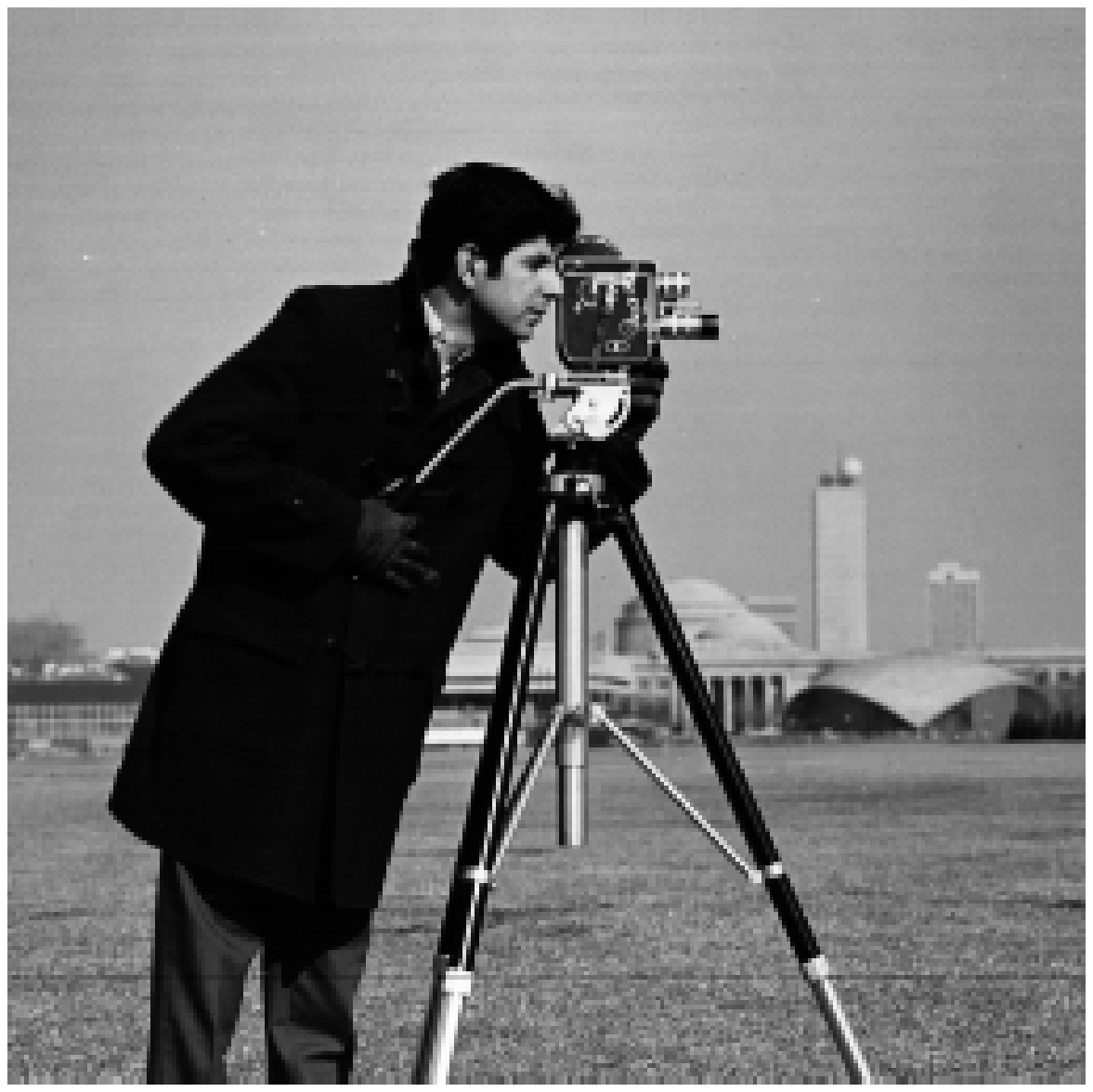}\\
\\
\includegraphics[width=.25\textwidth]{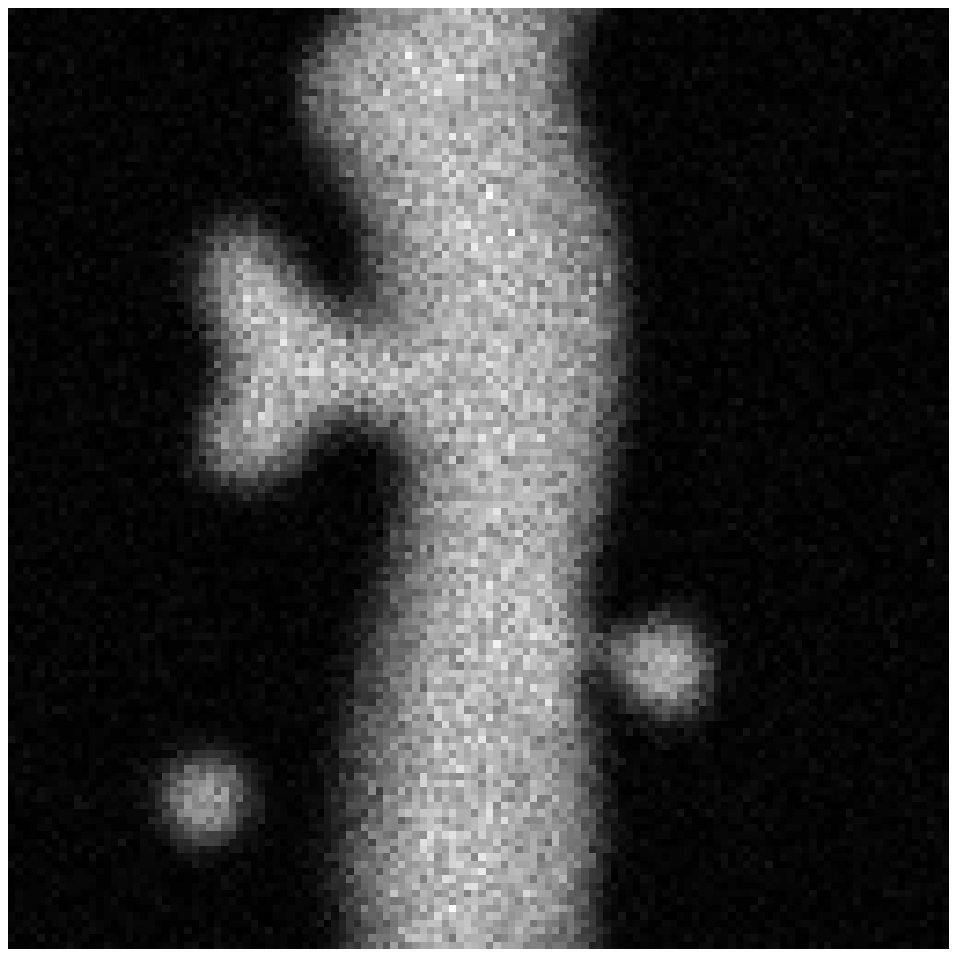}&
\includegraphics[width=.25\textwidth]{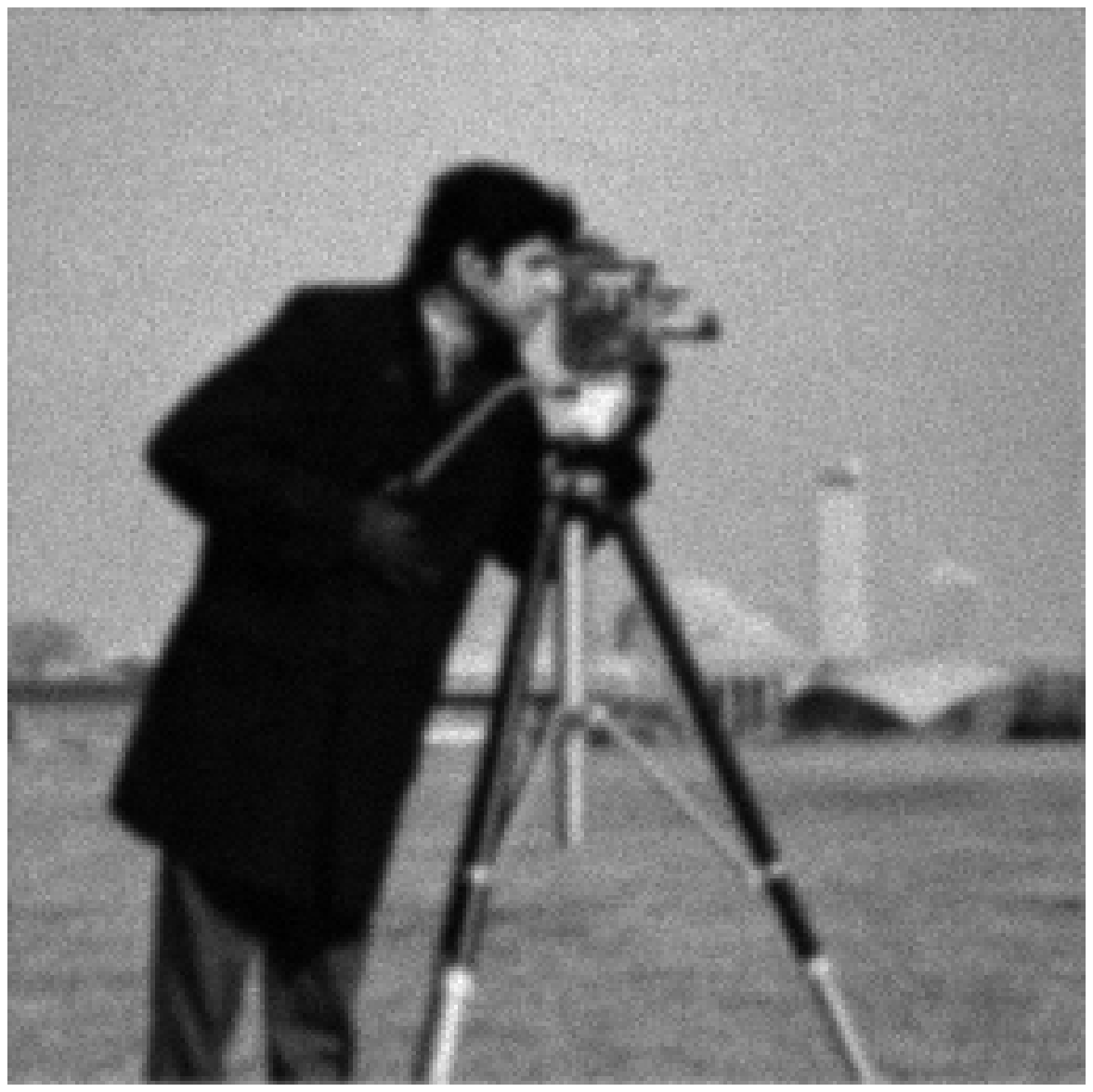}\\
Micro & Cameraman\\
 \end{tabular}
\caption{First row: original images for the two image deblurring test problems. Second row: blurred and noisy images for the two image deblurring test problems.}
\label{Immagini_KL+HS}
\end{center}
\end{figure}
\\
We compare the SFBEM approach with some other recent methods:
\begin{itemize}
 \item FISTA with backtracking \cite{Beck-Teboulle-2009b}, which can be considered as a special case of SFBEM where the scaling matrix is chosen at each iterate as the identity matrix. Actually our implementation slightly differs from the standard FISTA by the presence of the projection after the extrapolation step, which is needed when solving \eqref{min_KL+HS} since $\dom(f)$ does not coincide with $\R^n$;
 \item the scaled gradient projection method (SGP) \cite{Bonettini-etal-2009,Zanella-etal-2009}, which is a well known algorithm for differentiable constrained optimization problems. The SGP iteration has the form \eqref{SFB}, where the proximity operator reduces to the projection operator onto the constraints set. {Here the selection of steplength parameter $\alpha_k$ is based on the adaptive alternation of the Barzilai-Borwein rules proposed in \cite{Bonettini-etal-2009} and the value of $\lambda_k$ is computed by a line search procedure.} We point out that there exists a more recent variant of SGP \cite{Porta-Prato-Zanni-2015} which employs a different steplength selection rule. In order to avoid redundant results, we prefer to consider only the standard SGP approach as a comparative tool;
 \item the nonscaled version of SGP, hereafter indicated by GP.
\end{itemize}
The scaling matrix for SFBEM and SGP has been selected by exploiting the split gradient idea suggested in \cite{Lanteri-etal-2001} and based on a decomposition of the gradient into a {nonnegative} part and a negative one. {In particular, in \cite{Zanella-etal-2009} the authors show that the gradient of $f(x) = \KL(x)+\rho\HS(x)$ can be decomposed in the form
$$
-\nabla f(x) = U_{\KL}(x) + \rho U_{\HS}(x) - V_{\KL}(x) - \rho V_{\HS}(x)
$$
with $U_{\KL}, U_{\HS} \geq 0$ and $V_{\KL}, V_{\HS} >0$}. Then a possible scaling matrix is given by
\begin{equation}\label{Dk}
\ve{D}_k = {\rm{diag}}\left(\max\left(\frac{1}{\gamma_k},\min\left(\gamma_k,\frac{\ve{w}^{(k)}}{V_{\KL}(\ve{w}^{(k)}) + \rho V_{\HS}(\ve{w}^{(k)})}\right)\right)\right)^{-1}
\end{equation}
where the quotient is componentwise and $\ve w^{(k)}$ is equal to the previous iterate $x^{(k)}$ for SGP and to the extrapolated point $\y^{(k)}$ for SFBEM. Moreover we set $\displaystyle\gamma_k = \sqrt{1+\frac{10^{13}}{(k+1)^p}}, \ p=2.1$.

We remark that in this case the projection at Step 1 of SFBEM is independent on the scaling matrix $D_k$, since $Y$ is the nonnegative orthant, i.e. $ P_{Y,D} \equiv P_Y$. Thus we are allowed to first compute $\yk=P_Y(\xk+\beta_k(\xk-\xkb))$ then choosing $D_k$ depending on $\yk$ for updating $\xkk$.

Finally, the sequence $\{\beta_k\}$ employed in the definition of the extrapolation step for both FISTA and SFBEM has been chosen as in \eqref{thetak2} with $a=2.1$ in order to ensure the convergence of the sequence of the iterates.

The performance of the algorithms has been compared by evaluating their ability in reducing the objective function: in particular we computed an approximate solution $\xs$ of \eqref{min_KL+HS} by performing 20000 SGP iterations and, for any method, we consider at each iterate the relative difference between the objective function and the minimum value
\begin{equation}
\label{F_F_star}
\frac{F(x^{(k)}) - F(\xs)}{F(\xs)}
\end{equation}
Table \ref{Tab1_KL+HS} and Table \ref{Tab2_KL+HS} report the number of iterations and the computational time needed by each method to reduce the relative difference \eqref{F_F_star} below a certain tolerance $tol$. Since in both test problems the solution $\xs$ is unique, the relative minimization error
$ {\rm RME}(x^{(k)}) = \frac{\|x^{(k)} - \xs\|}{\|\xs\|}$ is also reported. The computational time presented is the average execution time (in seconds) over ten runs.

Figure \ref{F_F_star_KL+HS} shows the decrease of the relative differences \eqref{F_F_star} with respect to both the iteration number and the computational time.
We also observed that once the quantity \eqref{F_F_star} is below the tolerance $10^{-7}$,  all algorithms provide the same relative reconstruction error
${\rm RRE}(x^{(k)}) = \frac{\|x^{(k)} - x_{\rm true}\|}{\|x_{\rm true}\|}$, which measures the quality of the computed solution as an approximation of $x_{\rm true}$.
More precisely, this value is $0.088$ for the micro test problem $\bigl({\rm RRE}(b_{\rm \, micro}) = 0.195\bigr)$ and $0.087$ for the Cameraman one $\bigl({\rm RRE}(b_{\rm \, Cameraman}) = 0.121\bigr)$.
\begin{table}[htb!]
\begin{center}
\begin{tabular}{l|ccccccccc|}
                      & \multicolumn{9}{c|}{{\bf Micro}} \\
                      &\multicolumn{3}{c}{{\it tol} = $10^{-3}$} & \multicolumn{3}{c}{{\it tol} = $10^{-5}$} & \multicolumn{3}{c|}{{\it tol} = $10^{-7}$}\\

 & It. & RME & Time & It. & RME & Time & It.  & RME & Time \\
\hline
GP 		& 585 & 0.0414 & 5.69 & 2100 & 0.0077 & 21.31 & 3459 & 0.0014 & 34.70\\
SGP 		& 75  & 0.0261 & 0.72 & 203  & 0.0061 & 2.14  & 336  & 0.0016 & 4.48\\
FISTA 		& 223 & 0.0410 & 4.12 & 888  & 0.0048 & 15.26 & 3298 & 0.0005 & 56.05\\
SFBEM 		& 64  & 0.0202 & 0.84 & 188  & 0.0038 & 2.17  & 515  & 0.0004 & 6.72\\
\hline
\end{tabular}
\end{center}
\caption{Number of iterations and computational time required by each algorithm to reduce the relative difference \eqref{F_F_star} below given tolerances for the micro test problem. The corresponding RME and computational time (average over 10 runs) are also reported.}\label{Tab1_KL+HS}
\end{table}
\begin{table}[htb!]
\begin{center}
\begin{tabular}{l|ccccccccc|}
                      & \multicolumn{9}{c|}{{\bf Cameraman}} \\
                      &\multicolumn{3}{c}{{\it tol} = $10^{-3}$} & \multicolumn{3}{c}{{\it tol} = $10^{-5}$} & \multicolumn{3}{c|}{{\it tol} = $10^{-7}$}\\

 & It. & RME & Time & It. & RME & Time & It. & RME & Time \\
\hline
GP 		 & 1730 & 0.0134 & 76.76 & 4046 & 0.0021 & 164.07 & 5637 & 0.0003 & 229.16\\
SGP 		 & 241  & 0.0102 & 9.37  & 1178 & 0.0014 & 47.48  & 1671 & 0.0001 & 76.24\\
FISTA 		 & 226  & 0.0105 & 13.13 & 858  & 0.0011 & 54.46  & 3332 & 0.0001 & 220.71\\
SFBEM 		 & 42   & 0.0103 & 2.90  & 163  & 0.0012 & 10.31  & 705  & 0.0001 & 48.41\\

\hline
\end{tabular}
\end{center}
\caption{Number of iterations and computational time required by each algorithm to reduce the relative difference \eqref{F_F_star} below given thresholds for the Cameraman test problem. The corresponding RME and computational time (average over 10 runs) are also reported.}\label{Tab2_KL+HS}
\end{table}
\begin{figure}[htb!]
\begin{center}
\begin{tabular}{cc}
\includegraphics[width=.45\textwidth]{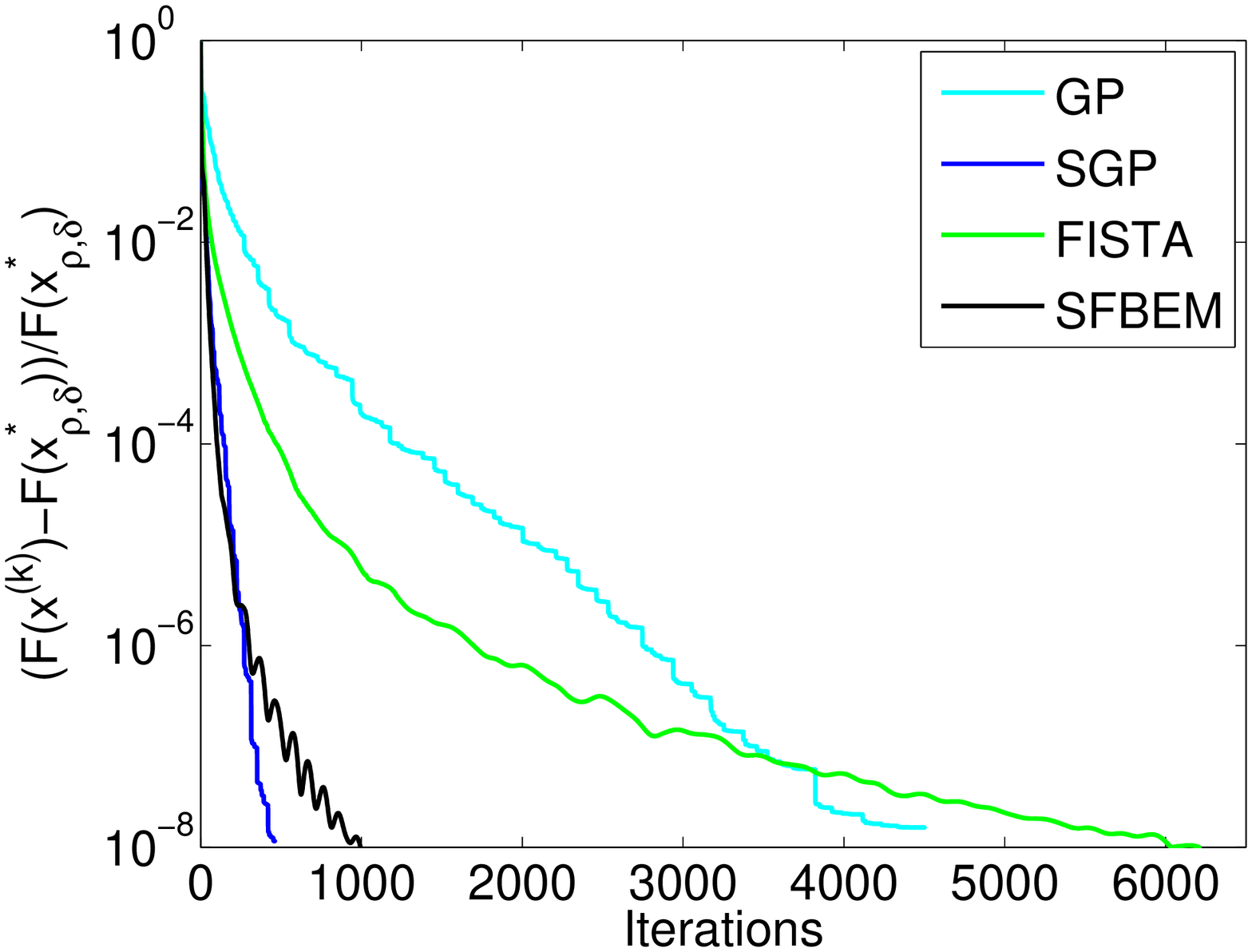}&
\includegraphics[width=.45\textwidth]{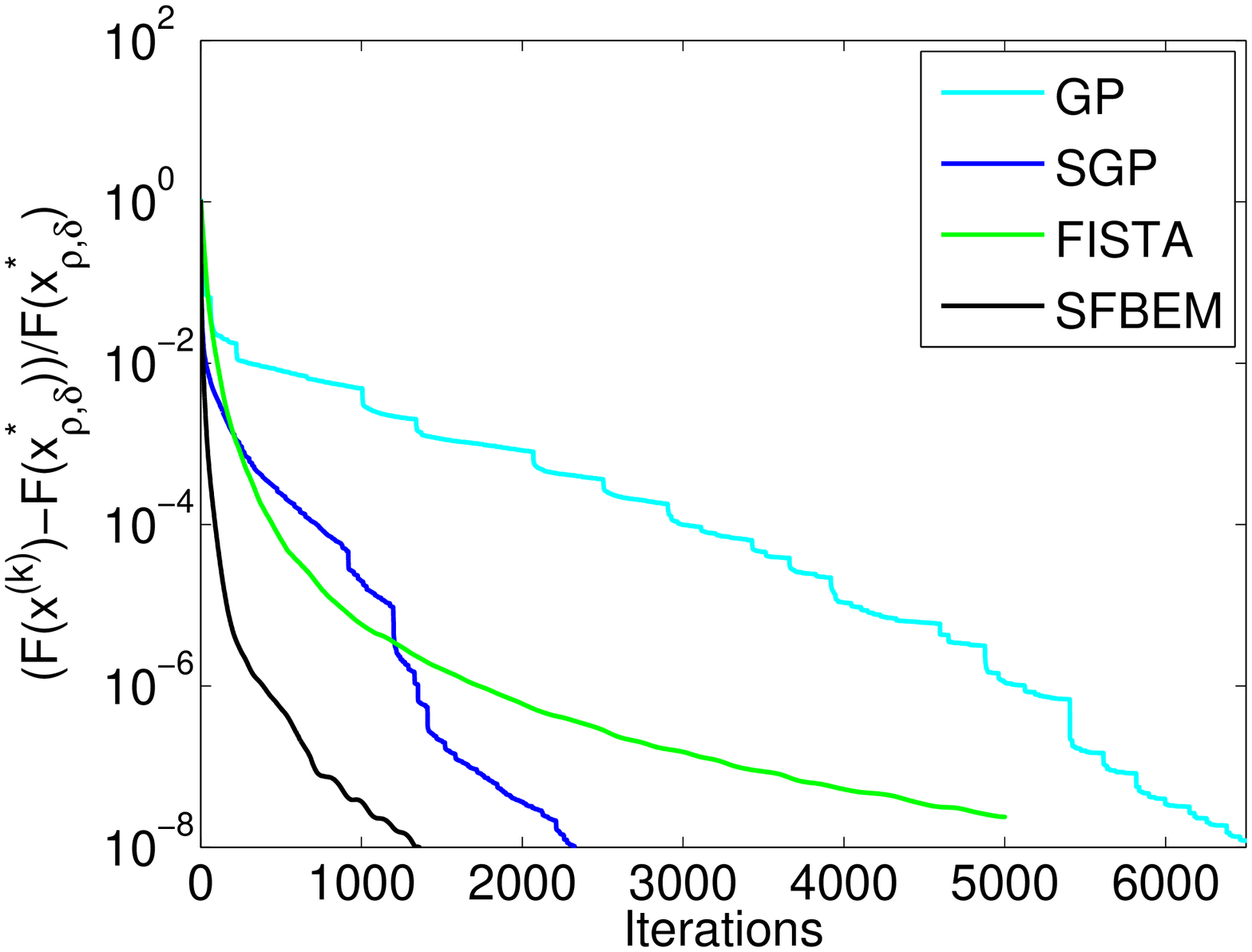}\\
\includegraphics[width=.45\textwidth]{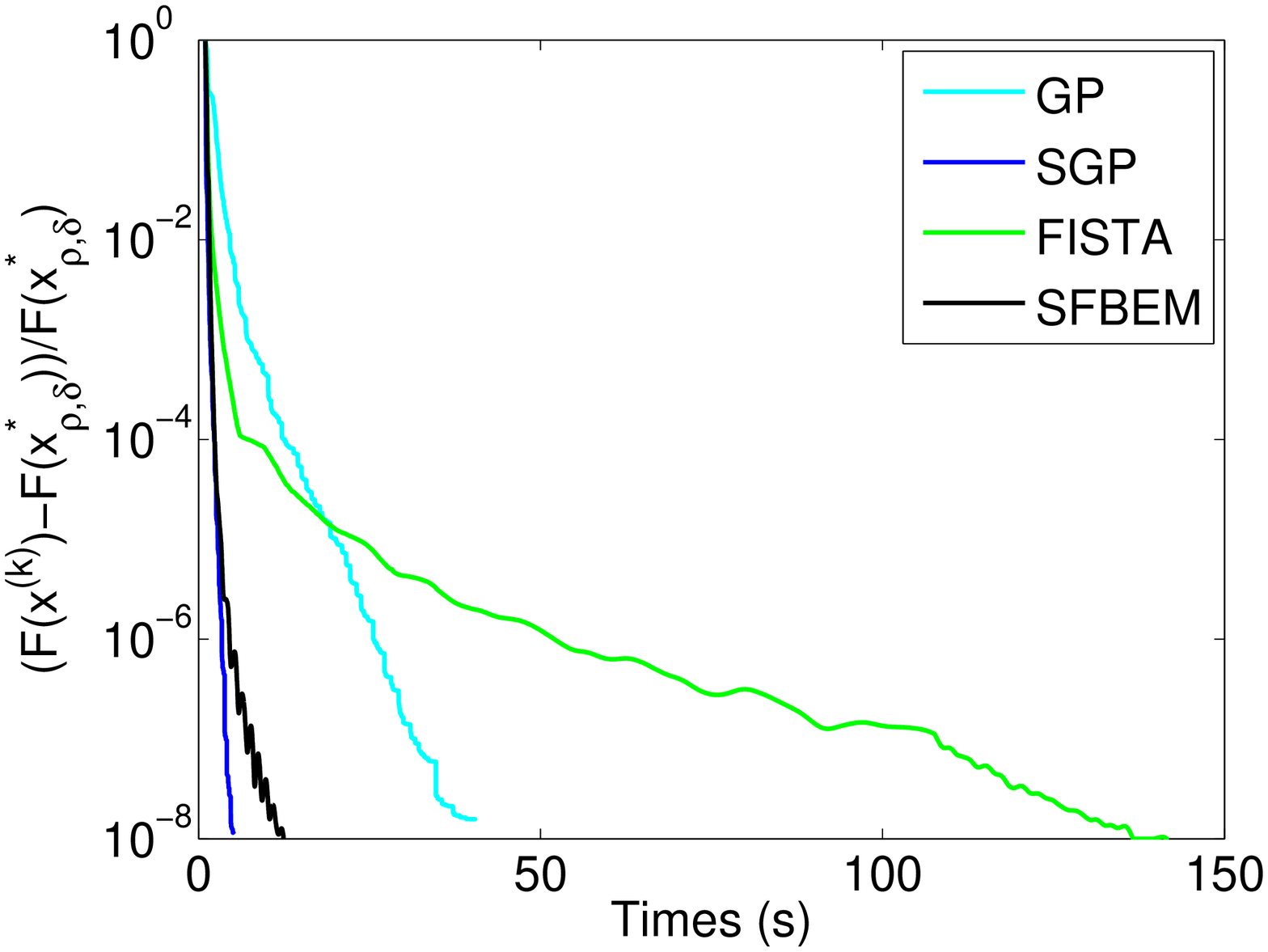}&
\includegraphics[width=.45\textwidth]{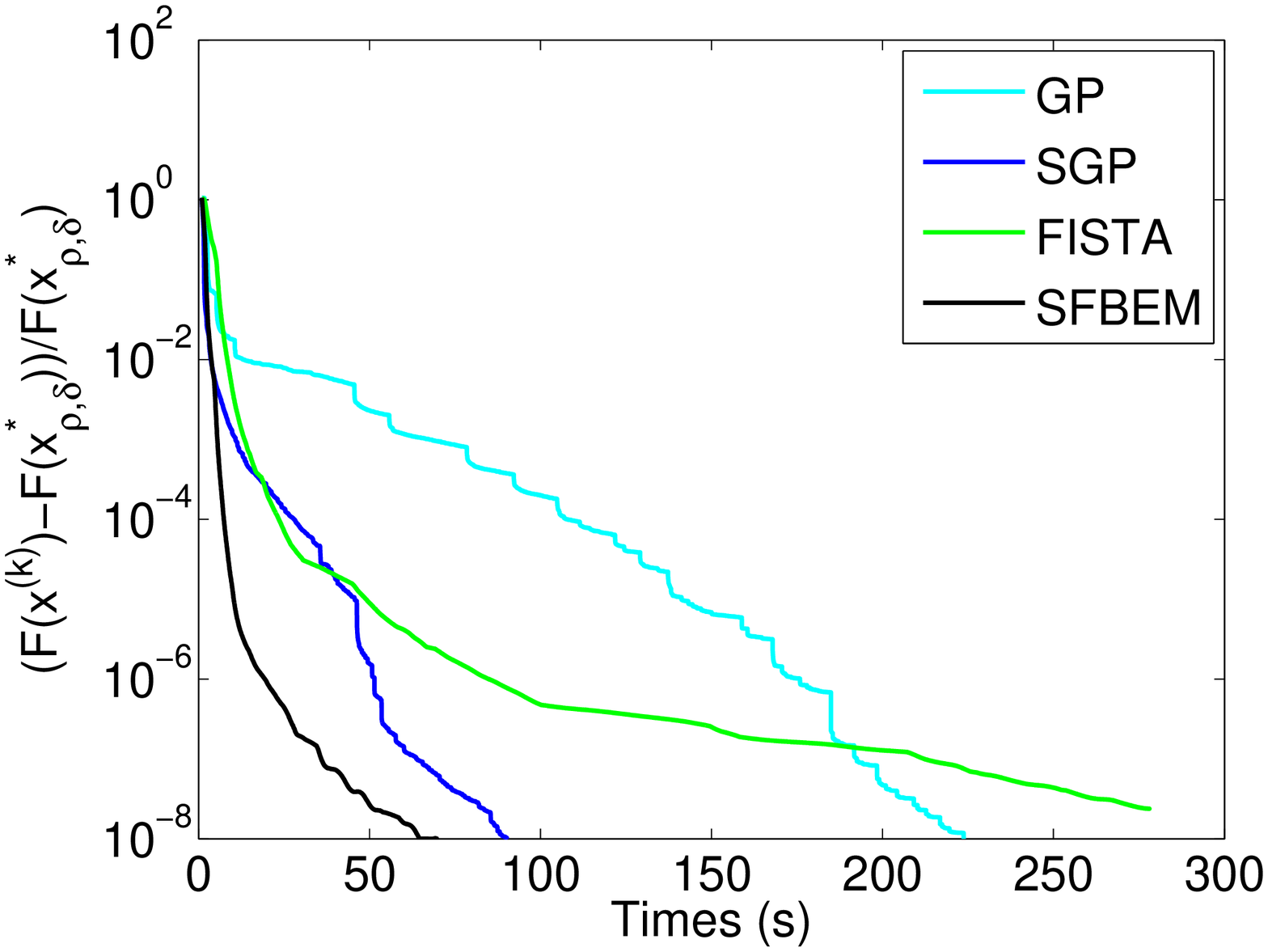}\\
Micro & Cameraman
\end{tabular}
\caption{Plots of the relative difference \eqref{F_F_star} with respect to the iterations number (top) and computational time (bottom) for the micro (left) and Cameraman (right) test problems.
}
\label{F_F_star_KL+HS}
\end{center}
\end{figure}

From the numerical results {shown} 
in Tables \ref{Tab1_KL+HS} and \ref{Tab2_KL+HS} and in Figure \ref{F_F_star_KL+HS}, it is possible to conclude that the presence of a non trivial scaling matrix makes the performances of SFBEM always superior to those of the nonscaled FISTA approach in terms of both number of iterations and computational time, while providing the same RRE and RME. Moreover, the comparison with SGP also supports the effectiveness of SFBEM: indeed, the performances of SFBEM are as good as those provided by SGP which is known in the literature as one of the most competitive algorithms to deal with image deblurring problems.
\subsection{Compressed sensing with Poisson noise}
As a second benchmark framework, we consider a compressed sensing problem which consists in recovering a sparse vector of {nonnegative} values starting from noisy measurements. More in detail, we assume that the observed data $b\in\mathbb{R}^m$ is the realization of a Poisson random variable with expected value given by $Ax_{\rm true} + bg$, where $x_{\rm true}\in\mathbb{R}^n$ is the signal of interest, $A\in\mathbb{R}^{m\times n}$ is the measurement matrix and $bg$ is a known background. As suggested in \cite{Raginsky10}, the true signal can be reconstructed by addressing a minimization problem of the form
\begin{equation}
\label{prob_compr_sen}
\min_{x\in\mathbb{R}^n} \KL(x) + \rho\|x\|_1 + \iota_{x\geq 0}(x)
\end{equation}
where $\KL$ is the {generalized} Kullback-Leibler divergence \eqref{KL_divergence}, the $\ell_1$ norm induces sparsity on the solution, $\rho$ is the positive regularization parameter and $\iota_{\ve{x}\geq \ve{0}}$ is the indicator function of the {nonnegative} orthant.

In this case, we set $f(\x) = \KL(\x) $ and $g(\x) = \rho\|\x\|_1 + \iota_{\x\geq 0}(\x)$ and $Y$ is the nonnegative orthant.
The operator $\projj\x$ associated to $g(x)$ can be computed in closed form \cite[Section II]{Harmany12}.

The numerical experiments are carried out on a test problem which has been generated with the following steps:
\begin{itemize}
 \item[(i)] a matrix $A\in\mathbb{R}^{1000\times5000}$ has been generated as {detailed} in \cite{Raginsky10} {so} that $A$ preserves both the positivity and the flux of any signal (i.e. if $z\geq0$ then $Az\geq0$ and $\sum_{i=1}^m (Az)_i \leq \sum_{i=1}^n z_i$);
 \item[(ii)] the signal to recover $x_{\rm true}\in\mathbb{R}^{5000}$ has all zeros except for 20 non-zero entries drawn uniformly in the interval $[0,10^5]$;
 \item[(iii)] the observed signal $b\in\mathbb{R}^{1000}$ has been obtained by corrupting the vector $Ax_{\rm true} + bg$ ($bg = 10^{-10}$) by means of the Matlab \verb+imnoise+ function.
\end{itemize}
The regularization parameter $\rho$ has been fixed equal to $10^{-3}$.

We compare SFBEM, FISTA with backtracking and the SPIRAL method developed in \cite{Harmany12} and designed to solve problems of the type \eqref{prob_compr_sen}. The SPIRAL approach {is a
forward-backward algorithm that employs a steplength selection strategy based on the Barzilai--Borwein rules \cite{Barzilai-Borwein-1988}; the convergence is guaranteed 
by a proper linesearch on the values of the objective function.} The Matlab code of SPIRAL is available on-line \cite{Spiral2012}. The scaling matrix for SFBEM has been selected by exploiting the already mentioned decomposition idea of the gradient of the differentiable part of the objective function: in particular, by writing the gradient of the KL functional as
$$
-\nabla \KL(x) = U_{\KL}(x) - V_{\KL}(x)
$$
with $U_{\KL}\geq 0$ and $V_{\KL}>0$, $D_k$ is defined as
$$
{D}_k = {\rm{diag}}\left(\max\left(\frac{1}{\gamma_k},\min\left(\gamma_k,\frac{\ve{y}^{(k)}}{V_{\KL}(\ve{y}^{(k)})}\right)\right)\right)^{-1}
$$
with $\displaystyle\gamma_k = \sqrt{1+\frac{10^{6}}{(k+1)^p}}, \ p=2.1$. The sequence $\{\beta_k\}$ used to update the extrapolation point has been chosen as in \eqref{thetak2} with $a=10$. The considered methods have been stopped when the relative distance between two successive iterations is less than $10^{-7}$. Table \ref{Tab_KL+L1} shows the performance of FISTA, SFBEM and SPIRAL in solving problem \eqref{prob_compr_sen} in terms of number of iterations and computational time (average over ten runs) to make the relative distance \eqref{F_F_star} smaller than prefixed thresholds. {For this test problem, }the minimum point $\xs$ has been computed by the SPIRAL method in 10000 iterations. Moreover, we report the relative reconstruction error and the relative minimization error. In order to better appreciate the results, Figure \ref{F_F_star_KL+L1} depicts the decreasing behavior of the objective function with respect to both the number of iterations
and the computational time. All the considered algorithms yield to the same value of the RRE equal to $0.075$.
\begin{table}[htb!]
\begin{center}
\begin{tabular}{l|ccccccccc|}
                      & \multicolumn{9}{c|}{{\bf Compressed sensing problem}} \\
                      &\multicolumn{3}{c}{{\it tol} = $10^{-3}$} & \multicolumn{3}{c}{{\it tol} = $10^{-5}$} & \multicolumn{3}{c|}{{\it tol} = $10^{-7}$}\\

 & It. & RME & Time & It. & RME & Time & It. & RME & Time \\
\hline
FISTA 	& 637  &  0.0111 &  15.75 & 933  &  0.0011 &  23.62 & 1412 &  0.0001 &  35.87\\
SFBEM 	& 369  &  0.0096 &  8.77  & 568  &  0.0011 &  14.37 & 806  &  0.0001 &  20.04\\
SPIRAL 	& 1180 &  0.0120 &  15.87 & 1309 &  0.0013 &  17.33 & 1379 &  0.0001 &  18.10\\
\hline
\end{tabular}
\end{center}
\caption{Number of iterations and computational time required by each algorithm to reduce the relative difference \eqref{F_F_star} below given tolerances for the compressed sensing test problem. The corresponding RME and computational time (average over 10 runs) are also reported.}\label{Tab_KL+L1}
\end{table}
\begin{figure}[htb!]
\begin{center}
\begin{tabular}{cc}
\includegraphics[width=.45\textwidth]{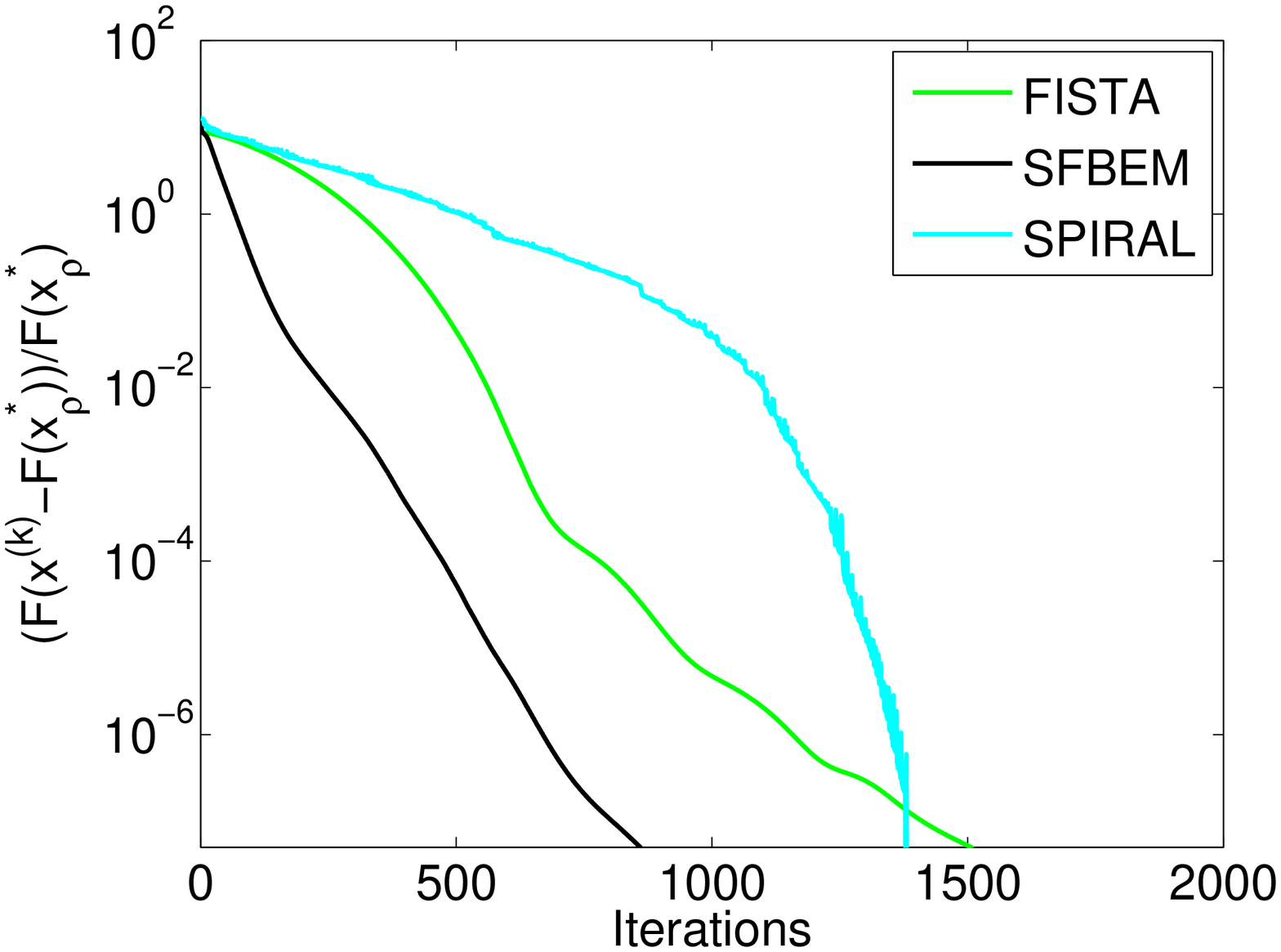}&
\includegraphics[width=.45\textwidth]{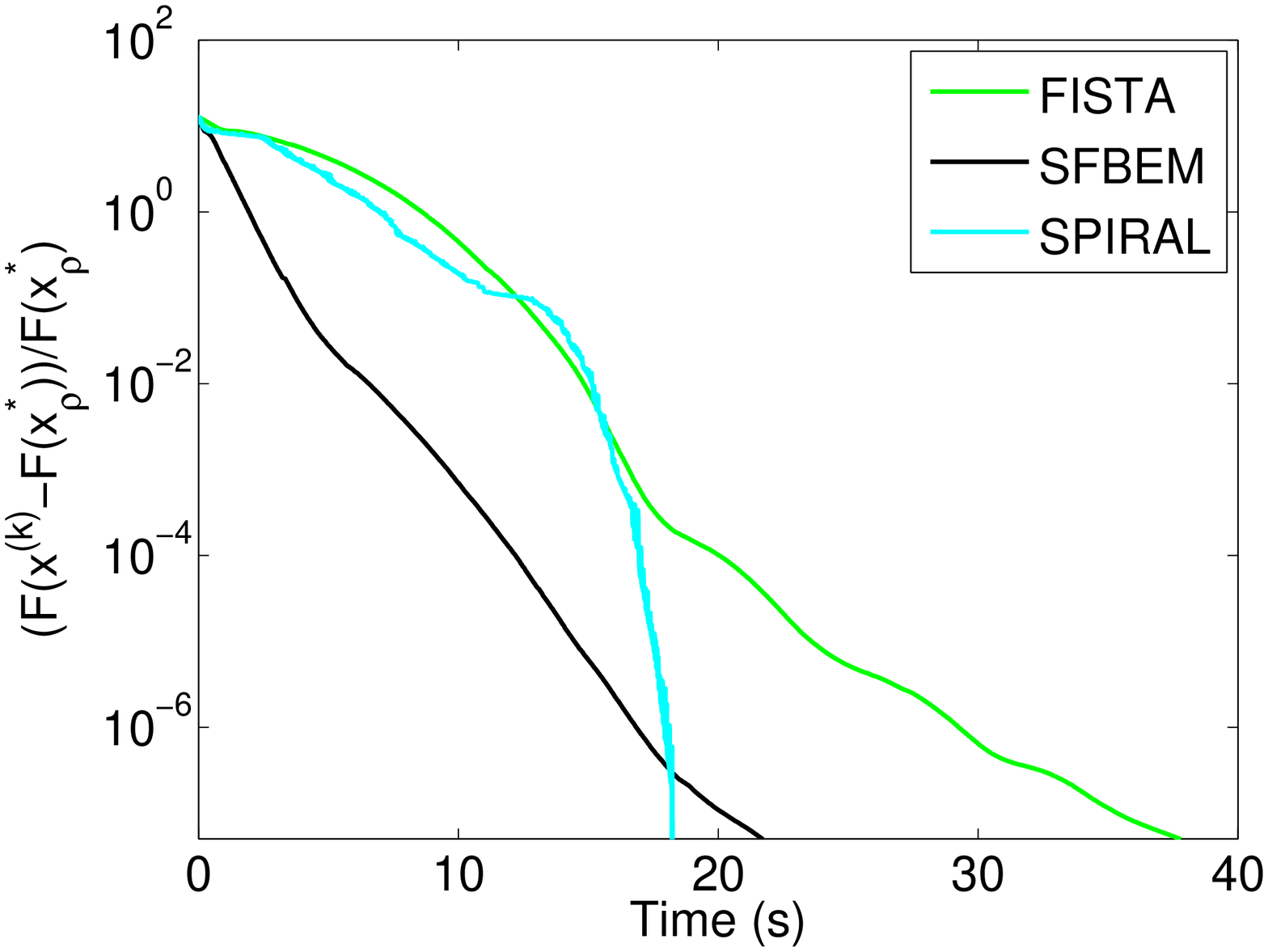}\\
\end{tabular}
\caption{Compressed sensing problem: relative difference \eqref{F_F_star} with respect to the iterations number (left) and computational time (right).
}
\label{F_F_star_KL+L1}
\end{center}
\end{figure}

The results obtained on the compressed sensing problem confirm the same conclusions in the image deblurring framework: the benefit of applying the SFBEM instead of FISTA is evident from the significant reduction of the number of iterations and  computational time, as reported in Table \ref{Tab_KL+L1}.
\subsection{Probability density estimation}
The last optimization problem we considered concerns the estimation of an unknown Gaussian mixture probability density. More in detail, if the sample $\{\tau_1, \tau_2, ..., \tau_n \, | \, \tau_i \in \mathbb{R}\}$ has been drawn from an unknown probability density function $\mu(\tau)$ which can be expressed as a Gaussian mixture, then a possible estimator has the form \cite{Girolami2003}
\begin{equation}
\label{estimator}
\hat{\mu}(\tau) = \sum_{i=1}^n x_i \kappa_{\sigma}(\tau,\tau_i)\,
\end{equation}
where $\kappa_{\sigma}(\cdot,\tau_i)$ is a Gaussian kernel with variance $\sigma$ and center $\tau_i$ and $x_i$ {is a suitable} coefficient. In \cite{Girolami2003} the authors proved that the weight vector $\x$ can be computed as a solution of the following minimization problem
\begin{equation}
 \label{prob_density_estimation}
 \min_{\x \in \mathbb{R}^n} \ \frac{1}{2} \x^TC\x - p^T\x + \iota_{\Delta_1^+}(\x)
\end{equation}
where
\begin{itemize}
 \item[(a)] the element $C_{i,j}$ of the matrix $C\in\mathbb{R}^{n\times n}$ is the Gaussian kernel of variance $2\sigma$, i.e. $C_{i,j} = \kappa_{2\sigma}(\tau_i,\tau_j)$;
 \item[(b)] the $i$-th component of the vector $p\in\mathbb{R}^n$ is defined as $\displaystyle p_i = \frac{1}{n}\sum_{j=1}^n \kappa_{\sigma}(\tau_i,\tau_j)$;
 \item[(c)] $\iota_{\Delta_1^+}$ is the indicator function of the simplex $\Delta_1^+ = \{\x\in\mathbb{R}^n \, | \, x_i\geq 0 \ \forall i, \ \ \sum_{i=1}^nx_i = 1\}$.
\end{itemize}
In this case we set $f(\x)= \frac{1}{2} \x^TC\x - p^T\x$, $g(\x)=\iota_{\Delta_1^+}(\x)$ and $Y = \R^n$. Thus, the operator $\projj\x$ consists in the projection onto the simplex $\Delta_1^+$. Such projection can be formulated as a root-finding problem and effectively computed by the secant-like algorithm proposed in \cite{Dai06}.

For the numerical experiments we analyzed the following Gaussian mixture
$$
\mu(\tau) = \frac{1}{5}\sum_{i=1}^5 \kappa_{\sigma_i}(\tau,c_i)
$$
with $\displaystyle\sigma_i = \sqrt[4]{\left(\frac{7}{9}\right)^{i-1}}$ and $\displaystyle c_i = 14\left(\left(\frac{7}{9}\right)^{i-1} - 1\right)$, $i=1,...,5$. The matrix $C$ and the vector $p$ have been generated with a sample of 1000 points drawn from $\mu$ by using the \verb+gmdistribution+ function of the Matlab Statistics and Machine Learning Toolbox.

The effectiveness of SFBEM has been evaluated in a comparison with FISTA with backtracking, SGP and GP.

The scaling matrix for SFBEM and SGP has been selected in the form
$$
\ve{D}_k = {\rm{diag}}\left(\max\left(\frac{1}{\gamma_k},\min\left(\gamma_k,\frac{\ve{w}^{(k)}}{C\ve{w}^{(k)}}\right)\right)\right)^{-1}
$$
with $\displaystyle\gamma_k = \sqrt{1+\frac{10^{10}}{(k+1)^p}}, \ p=2.1$ and $w^{(k)}$ equal to $y^{(k)}$ for SFBEM and $\xk$ for SGP. This choice of the scaling matrix mimics the split gradient based scaling proposed in \cite{Benvenuto-etal-2010} for quadratic problems. The extrapolation parameters sequence $\{\beta_k\}$ has been chosen as in \eqref{thetak2} with $a=2.1$. Table \ref{Tab_LS} reports the number of iterations and the computational time (an average value over ten runs) needed by the four methods to ensure a sufficient decrease of the distances \eqref{F_F_star}, where $\xs$ has been computed by means of 25000 FISTA iterations. GP and SGP do not succeed in satisfying the more restrictive threshold within the prefixed maximum number of iterations (25000). In Figure \ref{F_F_star_LS} we can appreciate the decrease of the objective function obtained by the considered algorithms and in Figure \ref{rec_LS} the reconstruction of the probability density function (pdf) through
the estimator \eqref{estimator} where
for simplicity we assume $\sigma = 1$.
\begin{table}[htb!]
\begin{center}
\begin{tabular}{l|cccccc|}
                      & \multicolumn{6}{c|}{{\bf Probability density estimation problem}} \\
                      &\multicolumn{2}{c}{{\it tol} = $10^{-3}$} & \multicolumn{2}{c}{{\it tol} = $10^{-5}$} & \multicolumn{2}{c|}{{\it tol} = $10^{-7}$}\\

 & It. &  Time & It. &  Time & It.  & Time \\
\hline
GP 		 & 111 & 0.51 & 20479 & 114.62 & - & - \\
SGP 		 & 90  & 0.98 & 4305  & 26.01  & - & - \\
FISTA 		 & 54  & 0.94 & 2141  & 16.54  & 21885 & 165.00 \\
SFBEM 		 & 53  & 0.67 & 810   & 6.46   & 3883  & 28.88\\
\hline
\end{tabular}
\end{center}
\caption{Number of iterations and computational time required by each algorithm to bring the relative difference \eqref{F_F_star} below given thresholds for the probability density estimation test problem.}\label{Tab_LS}
\end{table}
\begin{figure}[htb!]
\begin{center}
\begin{tabular}{cc}
\includegraphics[width=.45\textwidth]{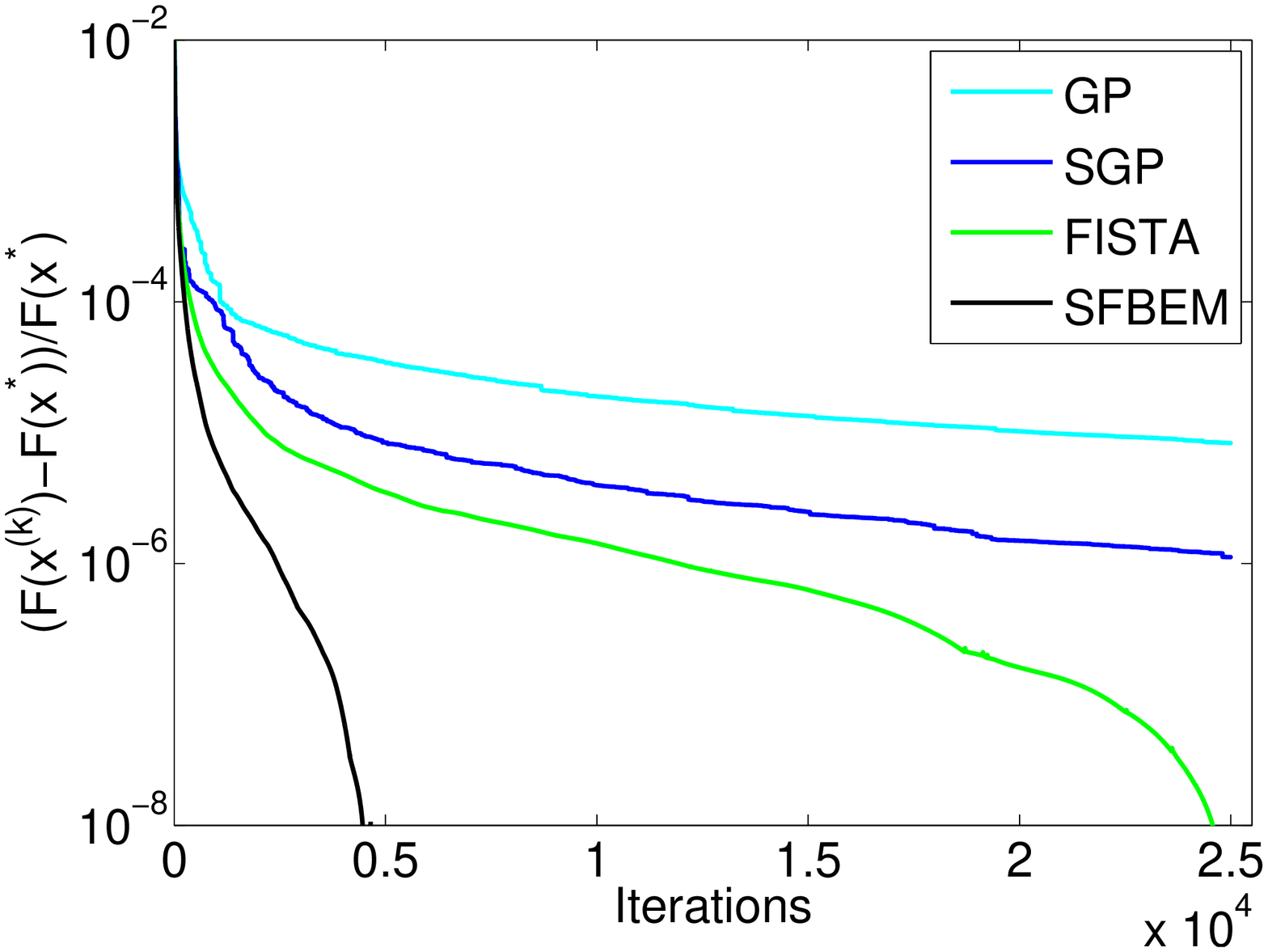}&
\includegraphics[width=.45\textwidth]{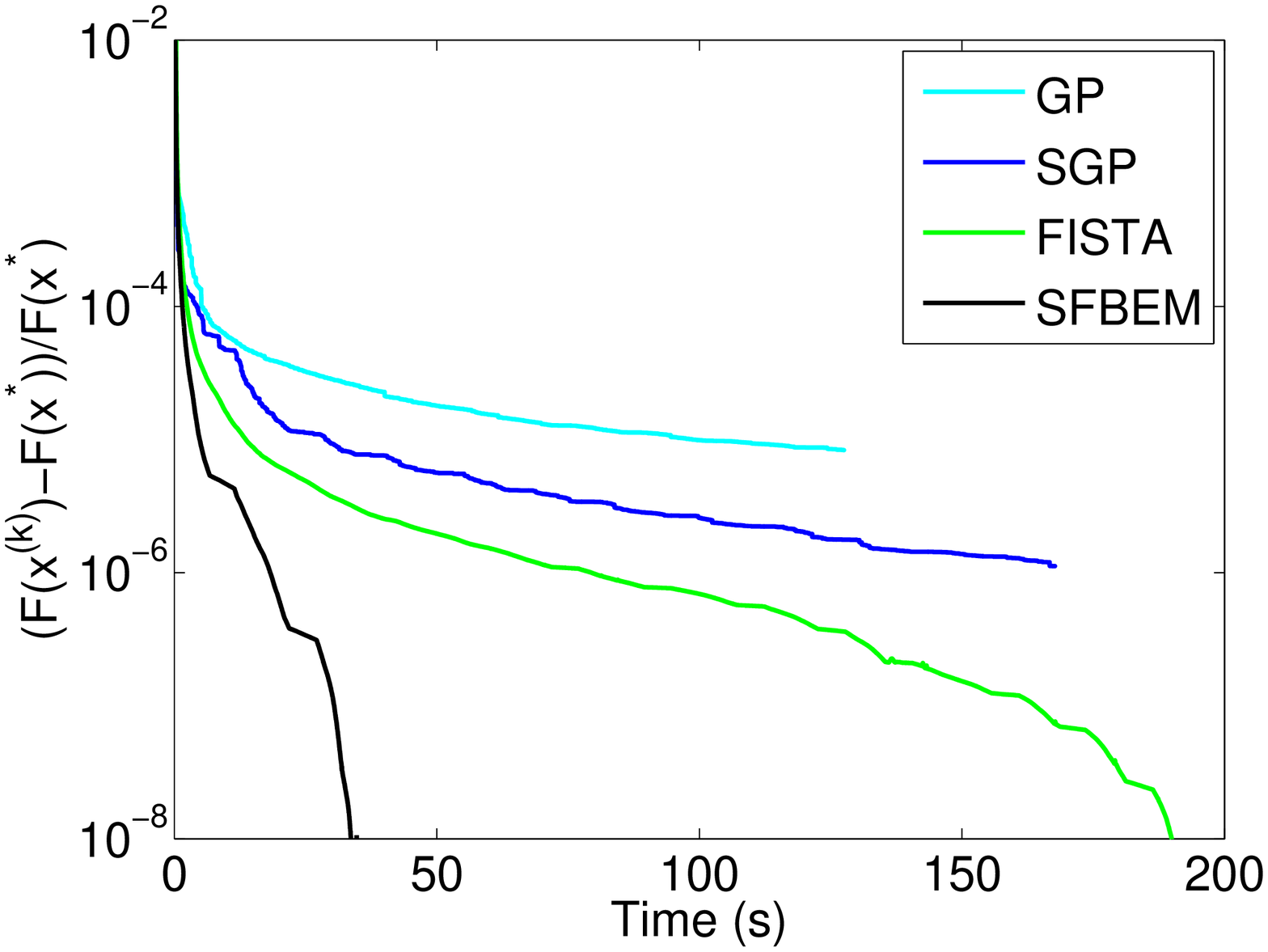}\\
\end{tabular}
\caption{Relative difference between the objective function values $F(\xk)$ provided by the different methods and the minimum value $F(\x^*)$.}
\label{F_F_star_LS}
\end{center}
\end{figure}
\begin{figure}[htb!]
 \begin{center}
 \includegraphics[width=.50\textwidth]{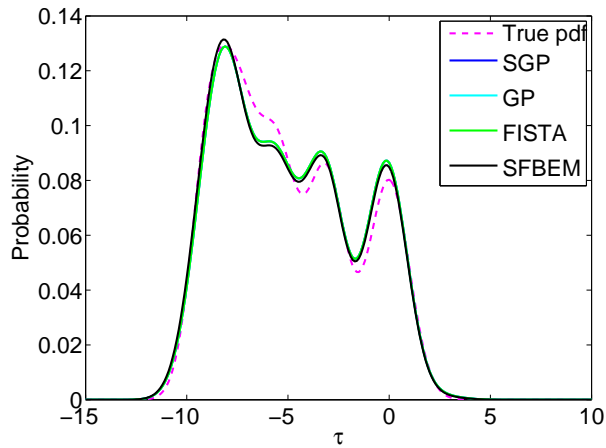}
  \caption{Density estimation results.}
\label{rec_LS}
 \end{center}
\end{figure}

The numerical experiments performed in the probability density estimation setting reinforce the validity of the proposed SFBEM scheme in comparison with the other state-of-the-art approaches we tested.
\section{Conclusions}\label{sec:concl}
In this paper we proposed a variable metric forward-backward method with extrapolation based on two fundamental ingredients:
a symmetric and positive definite scaling matrix multiplying the gradient of the differentiable part of the objective function, possibly capturing useful features of the problem to handle, and an inertial step which employs the information of the two last iterations in order to compute the new one.
A proper backtracking strategy ensuring a sufficient decrease of the objective function and suitable adaptive bounds on the scaling matrix allow to prove the convergence of the scheme to a minimizer of the considered problem. We also provided a convergence rate estimate which is similar to existing convergence rate results for nonscaled forward-backward algorithms with extrapolation.
Numerical experiments, carried out on optimization problems of different nature, showed very promising results in comparison with other algorithms which have already gained a great popularity in the literature. Future work will be addressed to analyze the possibility of introducing an inexact solution of the minimization problem which characterizes the backward step.
\bibliography{biblio_Silvia}
\bibliographystyle{plain}
\end{document}